\documentclass[11pt, a4paper]{article}
\usepackage{amssymb}
\usepackage{makeidx}
\usepackage{latexsym}
\newtheorem{theorem}{Theorem}
\newtheorem{axiom}{Axiom}
\newtheorem{metatheorem}[theorem]{Metatheorem}
\newtheorem{definition}[theorem]{Definition}
\newtheorem{metadefinition}[theorem]{Metadefinition}
\newtheorem{proposition}[theorem]{Proposition}

\newtheorem{lemma}[theorem]{Lemma}

\newtheorem{corollary}[theorem]{Corollary}

\newtheorem{ruleinf}{Rule of Inference}
\begin{document}

\def\P{\mathrm{P}}
\def\Field{\mathrm{Field}}
\def\First{\mathrm{First}}
\def\Last{\mathrm{Last}}
\def\Next{\mathrm{Next}}
\def\Prev{\mathrm{Prev}}
\def\Rank{\mathrm{Rank}}
\def\N{\mathcal{N}}
\def\M{\mathcal{M}}
\def\VN{\mathcal{VN}}
\def\Z{\mathcal{Z}}
\def\CH{\mathcal{CH}}
\def\ACK{\mathcal{ACK}}
\def\LEX{\mathcal{LEX}}
\def\Stage{\mathrm{Stage}}
\def\Lex{\mathrm{Lex}}
\def\<{\left \langle}
\def\>{\right \rangle}

\title{Natural Number Arithmetic in \\ the Theory of Finite Sets}
\author{J. P. Mayberry and Richard Pettigrew\footnote{Department of Philosophy, University of Bristol, 9 Woodland Road, Clifton, Bristol. BS8 1TB. \textsf{J.P.Mayberry@bris.ac.uk}, \textsf{Richard.Pettigrew@bris.ac.uk}}}
\date{}
\maketitle


\begin{abstract} We present here a finitary theory of finite sets and develop a theory of  `natural number arithmetic' that follows the example of Dedekind's theory of simply infinite systems expounded in \cite{rd}. In Dedekind's own theory, which rests on infinitary assumptions, all simply infinite systems are isomorphic, and can therefore all be regarded as concrete realisations of a single abstract system, \emph{the} natural number system. In our finitary version of Dedekind's theory, on the other hand,  simply infinite systems come in different lengths and have different closure properties with respect to arithmetical functions. Moreover, none of them contains enough numbers to specify the cardinalities of all finite sets. 
\end{abstract} 


The intuitive theory that underlies our formal theory is a finitary theory of finite sets. We shall call our theory (in both its formal and informal versions)  `Euclidean arithmetic' ($EA$), however,  because it represents a modern version of the approach to arithmetic taken by Euclid in Book VII of  the \emph{Elements}. In Euclid's arithmetic, a number (\emph{arithmos}) is a finite plurality composed of units; in other words, what we now call a finite set.\footnote{The abbreviation `EA' ( for `Elementary Arithmetic') is sometimes used to name the theory $I\Delta_0+\mathrm{exp}$. This threatens to cause confusion.  However, our notation is appropriate since there is a sense in which $I\Delta_0+\mathrm{exp}$ is the natural number arithmetic contained in our theory. This is explained in Section \ref{EAweakarith}, where it is also explained that there is another sense in which Euclidean Arithmetic is \emph{not} the set theory that corresponds to $I\Delta_0 + \mathrm{exp}$.}

From a foundational perspective we believe that Euclidean arithmetic is superior to its conventional rivals, which all, in their various ways, take the concept of natural number as a \emph{datum}, as something simply \emph{given} to our mathematical intuition. This is reflected in the fact that both mathematical induction and, at least, certain instances of recursion (in characterising addition and multiplication, for example) are built in as fundamental assumptions of the theory.

In Euclidean Arithmetic, by way of contrast, both induction and recursion are regarded as requiring justification on more basic, set-theoretical grounds. In taking this position, we follow the example of Dedekind \cite{rd} and indeed Frege \cite{gf}.  Perhaps our paper is best seen as carrying out Dedekind's project \cite{rd} in a finitary setting.

Euclidean Arithmetic is essentially a bounded quantifier version of the theory of hereditarily finite pure sets.\footnote{Thus all the sets in our theory are `built up' from the empty set in the familiar way. We could include urelements without affecting our principal results, but have opted not to do so in order to simplify our exposition.}  It is obtained from Zermelo-Fraenkel set theory by replacing the Axiom of Infinity with an Axiom of Dedekind Finiteness; including an Axiom of Transitive Closure; and restricting Subset Comprehension and Replacement to their bounded quantifier forms. 


In Sections \ref{formaltheory} - \ref{inducrec}, we present Euclidean arithmetic as a formal theory of finite sets, and develop as much of the theory as is indispensable to the presentation of our principal results. Sections \ref{linord} and \ref{inducrec} are expositions of the general theory of finite linear orderings: Section \ref{linord} introduces Kuratowski's definition of linear orderings, which is more convenient to our technical purposes than the usual one; in Section \ref{inducrec} we give a careful exposition of the theorems of induction and recursion for these orderings. In these sections we call attention to the importance of distinguishing between \emph{global} functions and relations, which are defined everywhere, and \emph{local} functions and relations, whose domains are finite sets.

In Section \ref{arithfuncrel},  we define the \emph{arithmetical} global functions and relations which, in effect, can be regarded as `operating' on cardinal equivalence classes, and thus correspond to the functions studied in conventional natural number arithmetic.

Section \ref{EAweakarith} surveys the relationship between $EA$ and weak fragments of Peano arithmetic.  We hope to dispel any impression that our results are essentially results in weak arithmetic, dressed up in unfamiliar language.

Sections \ref{sis} - \ref{wksubsysack}, which contain our principal results, deal with \emph{natural number systems} which are the analogues of Dedekind's \emph{simply infinite systems} treated in \cite{rd}. We show that in Euclidean Arithmetic these natural number systems come in differing, even incomparable, lengths, and satisfy different closure conditions with respect to familiar arithmetical functions. We define several hierarchies of natural number systems, whose constituent systems become longer and longer and are closed under more and more powerful arithmetical global functions.

We believe that these phenomena help to justify the foundational approach to arithmetic via the theory of finite sets that we have taken. In Euclidean Arithmetic, we are able to make important distinctions which cannot be made in more traditional foundational accounts that take the natural numbers as their starting point. This, of course, remains to be shown.

\section{The Formal Theory}\label{formaltheory}

The most fundamental formal version of Euclidean arithmetic is a free variable theory\ $EA_0$,  based on the following vocabulary:
\begin{enumerate}
\item[(i)] A constant, $\varnothing$, denoting the empty set.
\item[(ii)] First-order free variables, $a,b,c,\cdots$ (and with subscripts).
\item[(iii)] First-order binding variables (bound variables), for use in comprehension and replacement terms, $u,v,w,x,y,z, \cdots$ (and with subscripts).\footnote{We shall also sometimes use capital letters $S, T, U, V, etc.$ (and with subscripts) to stand for free variables or binding variables, as determined by the context.}
\item[(iv)] Propositional connectives, $\&$, $\lor$, $\lnot$, and $\supset$.
\item[(v)] Function symbols for pair set ($\{\_ \, ,\_\}$), power set ($\P$), union ($\bigcup$), and epsilon fan ($\epsilon$).\footnote{The \emph{epsilon fan}, $\epsilon(S)$, of set $S$ is the set of all linear orderings whose first term lies in $S$ and all of whose subsequent terms are members of their immediate predecessors.}
\item[(vi)] Term-forming operators for \emph{comprehension terms}
\[
\{x\in s:A[x/a]\}
\]
(where $A$ is an $EA_0$ formula and $s$ an $EA_0$ term) and \emph{replacement terms} 
\[
\{t[x/a]: x\in s\}
\]
(where $s$ and $t$ are $EA_0$ terms).
\item[(viii)] Relation symbols for membership ($\in$), identity ($=$),  and set inclusion ($\subseteq$).
\end{enumerate}
The syntactic formation rules for terms and formulas are the obvious ones. Note, however, that the presence of comprehension and replacement terms means that the notions of `term' and `formula' must be defined using a simultaneous induction.

We can think of an $EA_0$ term $t$ containing a single free variable $a$ as defining or determining a \emph{global function of one argument} which we may denote by the $\lambda$-expression $\lambda xt[x/a]$.\footnote{We can also form such expressions when the term $t$ has free variables in addition to $a$. These additional free variables we may think of as \emph{parameters} of the definition embodied in the $\lambda$-expression.} The function is \emph{global} since it is defined everywhere. Global functions of two or more arguments can be symbolised in the same way: $\lambda x_1\cdots x_\mathbf{n}t[x_1/a_1,\cdots,x_\mathbf{n}/a_\mathbf{n}]$ (where the $a_i$s are distinct free variables). Of course $\lambda xt[x/a](s)=t[s/a]$.

In the same way, we can think of an $EA_0$ formula  $A$ having a single free variable $a$ as determining a \emph{global relation  of one argument} which we may denote by $\lambda xA[x/a]$. Global relations of more than one argument can be symbolised in a similar way.


 We must emphasise that these $\lambda$ expressions are not part of the formal language but are merely auxilliary expressions. In any case, we shall rarely use them directly, for we introduce \emph{metavariables} for global functions and relations.
\begin{enumerate}
\item[(a)]Metavariables  for (global) functions of all degrees $>1$, $\sigma,\tau,\varphi,\cdots$ (and with subscripts).
\item[(b)]Metavariables  for (global) relations of all degrees $>1$, $\Phi,\Psi,\Upsilon,\cdots$ (and with subscripts). 
\end{enumerate}
Strictly speaking these metavariables do not range over global functions and relations, but rather over the $\lambda$-expressions that explicitly define them.    

It is an important feature of $EA_0$ that it permits us to define \emph{local quantifiers}:  
\[
(\forall x\in S)A[x/a]=_{\mathrm{df.}}S\subseteq\{x\in S:A[x/a]\}
\]
and
\[
(\exists x\in S)A[x/a]=_{\mathrm{df.}}\lnot(\forall x\in S)\lnot A[x/a]
\]
using the comprehension operator. Note that formulas containing these `quantifiers' are, in fact, quantifier-free in the usual sense. These definitions  permit us to avoid the inclusion of (global) quantifiers in our list of primitive symbols for $EA_0$.

The axioms and rules of $EA_0$ comprise those of the classical propositional calculus, the usual substitution rule for free variables,
together with the following set-theoretical axioms and rule of inference.  First, the axioms that characterize the primitive constants and functions of $EA_0$:
\begin{axiom}[Empty Set] $a \not \in \varnothing$ \end{axiom}
\begin{axiom}[Pair Set] $a \in \{S, T\} \equiv (a = S \vee a = T)$. \end{axiom}
\begin{axiom}[Power Set] $a \in \P(S) \equiv a \subseteq S$. \end{axiom}
\begin{axiom}[Union] $a \in \bigcup S \equiv (\exists X \in S)[a \in X]$. \end{axiom}
\begin{axiom}[Schema of Subset Selection] Suppose $\Phi$ is a global relation of one place.\footnote{This preliminary sentence is logically superfluous, since the axiom schema being laid down conforms to our conventions on the use of metavariables. A similar remark applies to the Axiom Schema of Replacement.}  Then
\[
a \in \{x \in S : \Phi(x)\} \equiv a \in S \ \&\  \Phi(a)
\]
\end{axiom}
\begin{axiom}[Schema of Replacement] Suppose $\varphi$ is a global function of one argument.  Then
\[
a \in \{\varphi(x) : x\in S\} \equiv (\exists x\in S)[a=\varphi(x)]
\]
\end{axiom}

Next, the axiom and rule of inference that together characterize the relation, $\subseteq$, of subset inclusion:
\begin{ruleinf}[$\subseteq$-introduction] Suppose $s$ and $t$ are terms of $EA_0$.  Then, from
\[
a\in s\supset a\in t
\]
as premise infer
\[
s\subseteq t
\]
as conclusion, provided that the free variable $a$ does not occur in the conclusion or in any hypothesis on which the premise depends.
\end{ruleinf}
\begin{axiom}[$\subseteq$-elimination] $S \subseteq T\ \&\  a \in S\, .\supset.\, a \in T$.
\end{axiom}
And, finally, the axioms that determine the nature of the sets of our theory.
\begin{axiom}[Extensionality] $a \subseteq b\ \&\  b \subseteq a\ . \supset.\, a = b$ \end{axiom}
There are two further axioms, but we postpone their statement until Section \ref{tworemax}.

With these axioms and rules in place, it is a straightforward matter to introduce all the familiar operations of set theory, such as Boolean union and intersection, set difference, ordered pair, Cartesian product, the concept of local relation and local function, and so on. Similarly, the global relations of cardinal equivalence ($\simeq_\mathrm{c} $) and cardinal inequality ($\leq_\mathrm{c} $) can be defined and their basic properties established. This allows us to state our final axiom, which replaces conventional, infinitary ZF's Axiom of Infinity.
\begin{axiom}[Dedekind Finiteness] $a\subsetneqq b \,.\supset.\,a<_\mathrm{c} b$.
\end{axiom}

Using the axiom of Dedekind finiteness, we can establish a fundamental theorem that is a form of induction for finite sets. The \emph{one point extension of $S$ by $a$} is just the set $S\cup\{a\}$. 

\begin{theorem}[Local One Point Extension Induction]\label{onepointext} Let $S$ be any set and $T$ any set of subsets of $S$ satisfying
\begin{enumerate}
\item[\emph{(i)}] $\varnothing \in T$.
\item[\emph{(ii)}] For every $X \in T$ and $y \in S$, $X\cup\{y\} \in T$.
\end{enumerate}
Then $S \in T$.
\end{theorem}

\noindent\emph{Proof}:  This theorem can be stated formally and proved formally in $EA_0$.  Similar remarks apply below.

Let $S$ and $T$ be sets satisfying the hypotheses of the Theorem, and suppose, by way of contradiction, that $S\not\in T$. Then for every $X\in T$ there must be a $y\in S$ such that $y\not\in X$. But, by condition (ii), $X\cup\{y\}\in T$.  Thus, for every $X\in T$ let $\hat{X}=\{Y\in T: Y\simeq_\mathrm{c}  X\}$ and let $\hat{T}=\{\hat{X}:X\in T\}$. Now set $f=\{(\hat{X},\widehat{X\cup\{y\}}):y\in (S-X),\,\,\hat{X}\in \hat{T}\}$. Then $f:\hat{T}\to\hat{T}$ and $f$ is injective. But $f$ is not surjective since $\widehat{\varnothing}\notin f ``\hat{T}$. This violates the \emph{axiom of Dedekind finiteness}.\hfill$\Box$ \medskip

The local version of one point extension induction yields the following theorem as an easy corollary.

\begin{theorem}[Global One Point Extension Induction]\label{glonepointext} Let $\Phi$ be a global relation of one argument. Then from the premises
\begin{enumerate}
\item[\emph{(i)}]$EA_0\vdash \Phi(\varnothing)$.
\item[\emph{(ii)}]$EA_0\vdash \Phi(S) \supset \Phi(S\cup\{b\})$, where the free variable $S$ does not occur in $\Phi$.
\end{enumerate}
we may infer that $EA_0\vdash \Phi(S)$.
\end{theorem}


Define a \emph{choice function} for a set $S$ to be a function $f: (\P(S)-\{\varnothing\})\to S$ such that $f`X\in X $, for all non-empty subsets, $X$, of $S$.

\begin{theorem}[The Principle of Choice]\label{axchoice} Every set has a choice function.
\end{theorem}  

\noindent\emph{Proof}: A straightforward application of One Point Extension Induction. \hfill $\Box$\medskip



Before we can complete the presentation of $EA_0$ we must first introduce the concept of linear ordering.

\section{Linear orderings}\label{linord}

We need to develop the theory of linear orderings in $EA_0$, not only in order to complete the presentation of the axioms, but also because linear orderings are central to the theory of natural number systems, which forms the core of our presentation.

Although we follow the standard treatment of local functions as sets of Kuratowski ordered pairs, our treatment of linear orderings deviates from the usual presentation:  indeed, we follow Kuratowski's own treatment in \cite{kk}.

The idea is this:  Suppose $a_i \neq a_j$ unless $i = j$; then the linear ordering 
\[
[a_1, a_2, \cdots, a_{n-1}, a_n]
\]
is defined to be the set of all (fields of) its non-empty initial segments
\[
 \{\{a_1\}, \{a_1, a_2\}, \cdots, \{a_1, a_2, \cdots, a_{n-1}\}, \{a_1, a_2, \cdots, a_{n-1}, a_n\}\}
\]
We call the set of terms of a linear ordering $L$ the \emph{field of $L$} and denote it $\Field(L)$:  thus, $\Field([a_1, \cdots, a_n]) = \{a_1, \cdots, a_n\}$.

Of course, this description is merely heuristic.  Now, we give the rigorous definition:
\begin{definition}[Linear Ordering]\index{linear ordering}
By definition, $L$ is a linear ordering, if
\begin{enumerate}
\item[\emph{(i)}] $L$ does not contain the empty set;
\item[\emph{(ii)}] If $L$ is non-empty, then $L$ contains a singleton;
\item[\emph{(iii)}] For every $x \in L$, either $x = \bigcup L$, or there is $y \in \bigcup L$ such that $y \not \in x$ and $x \cup \{y\} \in L$; and
\item[\emph{(iv)}] For all $x, y \in L$, either $x \subseteq y$ or $y \subseteq x$.
\end{enumerate}
\end{definition}
When linear ordering is defined in this way, the field of a linear ordering is simply the union of that linear ordering:  thus, by definition,
\[
\Field(L) =_{\mathrm{df.}}  \bigcup L.
\]
\begin{definition} If $L$ is a linear ordering, and $a, b \in \Field(L)$, then, by definition,
\[
a <_L b \mbox{ if, and only if, } (\exists x \in L)[a \in x\ \&\  b \not \in x]
\]
\end{definition}
Some remarks about these definitions:
\begin{enumerate}
\item[(1)]  Clause (i) rules out $\{\varnothing, \{a, b\}\}$; clause (ii) rules out $\{\{a, b\}, \{a, b, c\}\}$; clause (iii) rules out $\{\{a\}, \{a, b, c\}\}$; and clause (iv) rules out such sets as $\{\{a, b\}, \{c, d\}\}$ and $\{\{a, b\}, \{a, c\}\}$.
\item[(2)]  Given this definition, it follows from clause (iv) that the field of a linear ordering is also the largest set in that linear ordering.  Further, it is possible to prove, by One Point Extension Induction and clauses (i), (ii), and (iii) that the largest set in a linear ordering is of the same cardinality as the whole linear ordering.  Putting these together, if $L$ is a linear ordering, then
\[
\Field(L) \simeq_\mathrm{c}  L.
\]
(It is because of this feature of Kuratowski linear orderings that we have chosen to use them.) We say that a linear ordering $L$ is a linear ordering \emph{of} the set $S$ if, and only if, $S=\Field(L)$. Thus if $L$ is a linear ordering of $S$, $L\simeq_\mathrm{c}S$.
\item[(3)] Suppose we define the global relation $\lambda x, y\ x \simeq_\mathrm{o} y$ in the usual way:  that is, by definition, $L_1 \simeq_\mathrm{o} L_2$ if, and only if, there is a bijection $f : \Field(L_1) \to \Field(L_2)$ such that, for every $x$ and $y$ in $\Field(L_1)$, $x <_{L_1} y$ if, and only if, $f ` x <_{L_2} f ` y$.  Then, we can show that, for any two linear orderings, $L_1$, $L_2$,
\[
L_1 \simeq_\mathrm{o} L_2 \mbox{ if, and only if, } L_1 \simeq_\mathrm{c} L_2
\]
\item[(4)] By our definition $\varnothing$ is a linear ordering:  when we are considering it as a linear ordering, we call it the \emph{empty ordering} and denote it by `$[\, ]$'.
\item[(5)] Since in $EA_0$ any set, and thus any linear ordering, must be Dedekind finite, it is possible to define global functions
\[
L\mapsto \First(L);\ \ L\mapsto \Last(L); \ \ L,x\mapsto \Next_L(x);\ \ L,x\mapsto \Prev_L(x)
\]
which all have the properties that their names indicate and whose values all lie in $\Field(x)$.\footnote{$\Next_L(\Last(L))=\Last(L)$ and $\Prev_L(\First(L))=\First(L)$.}
\item[(6)] It is also possible to define the global relation
\[
\lambda xy[x\subseteq_*y]
\]
which means that $L_1$ and $L_2$ are linear orderings and $L_1$ is an initial segment of $L_2$.  In fact, it is a consequence of the Kuratowski definition of linear ordering that if $L_1$ and $L_2$ are both linear orderings, then  $L_1 \subseteq_* L_2$ if, and only if, $L_1 \subseteq L_2$.

Given a linear ordering, $L$, we define $\mathrm{InSeg}(L)$ to be the set of all \emph{proper} initial segments of $L$.
\end{enumerate}

\section{The two remaining axioms}\label{tworemax}

We return to the properties of linear orderings in Section \ref{inducrec}. But first we turn to the statement of the two remaining axioms of $EA_0$. We begin with the axiom of $EA_0$ that characterizes the \emph{epsilon fan} of a set. First, we need a definition.

\begin{definition} If $L$ is a linear ordering, we say that \emph{$L$ is an $\epsilon$-chain for $S$}  if
\begin{enumerate}
\item[\emph{(1)}] $L$ is non-empty;
\item[\emph{(2)}] $\First(L) \in S$; and
\item[\emph{(3)}] For all $x \in \Field(L)$ except $\Last(L)$, $\Next_L(x) \in x$.
\end{enumerate}
\end{definition}
We can now give a precise formulation of the Axiom of the Epsilon Fan:

\begin{axiom}[Epsilon Fan] $a\in\epsilon(b)\equiv a$ is an $\epsilon$-chain for $b$.
\end{axiom}

This allows us to define the \emph{transitive closure}, $TC(S)$ of a set $S$:
\begin{definition}$\mathrm{TC}(S)=_{\mathrm{df.}}\bigcup(\{\mathrm{Field}(x):x\in\epsilon(S)\})$
\end{definition}
The transitive closure of a set comprises all those sets that go, directly or indirectly, to make up the set.  With this in hand, we can now state the Axiom of Foundation:

\begin{axiom}[Foundation] $a\in \mathrm{TC}(b)\supset b\not\in \mathrm{TC}(a)$\end{axiom}
This axiom says that we can't have sets $S$ and $T$ each of which is used in `building up' the other.

This completes our presentation of the axioms and rules of inference of our quantifier-free theory of finite sets:  together, Axioms 1 to 11 and Rules of Inference 1 and 2 constitute $EA_0$.

$EA_0$ can be extended to a classical first-order theory $EA_1$ by adding global quantifiers $\forall$ and $\exists$, together with the usual axioms and rules. The definition of formula must be extended to accommodate the new global quantifiers, but the definition of term remains unchanged. In particular, \emph{the principal formula $A$ in the comprehension term $\{x\in s:A[x/a]\}$ must be a formula of $EA_0$ and thus free of the new unbounded quantifiers}. $EA_1$ can be shown to be conservative over $EA_0$.

\section{Induction and recursion}\label{inducrec}

In ordinary natural number arithmetic, both induction and recursion apply to the whole set of natural numbers, and induction at least is included among the basic axioms in formal treatments of arithmetic.\footnote{As noted above, two particularly important instances of recursion, namely, the recursion equations for addition and multiplication, are usually included among the formal axioms.}  Since we are approaching arithmetic via the theory of finite sets, and must prove our results using only set-theoretical principles, we shall follow Dedekind in regarding proof by induction and definition by recursion as propositions that we must prove, rather than as assumptions we may make in our axioms. Of course, since we are dealing only with \emph{finite} sets we cannot employ all of Dedekind's ideas.

The fundamental technical problems with induction and recursion occur already when we consider proof by induction and definition by recursion along a linear ordering, which, being a set, is necessarily finite in our case.

\begin{theorem}[Induction along a Linear Ordering]\label{locind} Suppose $\Phi$ is a global relation \emph{(}i.e. defined by a $EA_0$ formula\emph{)}. Then if
\begin{enumerate}
\item[\emph{(i)}]$L$ is a linear ordering,
\item[\emph{(ii)}]$\Phi(\First(L))$, and
\item[\emph{(iii)}] For all $x\in\mathrm{Field}(L)$ except $\Last(L)$, $\Phi(x) \mbox{ implies } \Phi(\Next_L(x))$,
\end{enumerate}
then, for all $x \in \mathrm{Field}(L)$, $\Phi(x)$.
\end{theorem}

\noindent\emph{Proof}:  As in the proof of induction in ZF, we assume the premises and suppose, for a contradiction, that there is $x \in \Field(L)$ such that $\neg \Phi(x)$. We let
\[
S = \{x \in \Field(L) : \neg \Phi(x)\}
\]
and derive a contradiction in the obvious way. \hfill$\Box$ \medskip


The proof of Theorem  \ref{locind} depends essentially on the fact that the induction formula belongs to $EA_0$.   If it were to include unbounded quantifiers, then this proof could not be given since it would be illegitimate to form the term
\[
S = \{x \in \Field(L) : \neg \Phi(x)\}
\]

\begin{theorem}[Local Recursion along a Linear Ordering]\label{locrec}Let $L$ be a linear ordering, $S$ a set, $a$ an element of $S$, and $g:S\to S$ a local function. Then there is exactly one local function $f:\mathrm{Field}(L)\to S$ such that
\begin{enumerate}
\item[\emph{(i)}]$f`\First(L)=a$ and
\item[\emph{(ii)}]For all $x\in\mathrm{Field}(L)$ except for $\Last(L)$, $f`\Next_L(x)=g` f`x$.
\end{enumerate}
\end{theorem}

\noindent\emph{Proof}: Again, the proof is similar to the proof of recursion in ZF:  we prove, by induction along $L$, that there exists a partial solution to the recursion equations restricted to the initial segment of $L$ up to that point.  Since we can specify in advance a set in which all such partial solutions will lie, the induction formula is bounded and thus $EA_0$:  it is the set of all local functions from the field of an initial segment of $L$ into $S$.  Uniqueness is established by a second bounded induction. \hfill $\Box$\medskip


It is clear that such a proof would not work for the global version of this theorem in which the local function $g:S\to S$ is replaced by a global function $\varphi$. In the global case, we cannot specify in advance a set that contains all partial solutions to the recursion equations:  prior to carrying out the recursion, we cannot specify a set $S$ such that all the partial solutions belong to the set of local functions from initial segments of $L$ into $S$.

The failure of induction along a linear ordering with respect to formulas containing global quantifiers and the failure of recursion along a linear ordering with respect to an arbitrary global function play a crucial part in the remainder of our exposition. These failures are caused by the restriction of the formula $A$ in comprehension terms of the form $\{x\in S:A[x/a]\}$, which is intended to reflect Zermelo's requirement that such formulas be \emph{definite}. Such considerations presumably motivate the restriction of the induction axioms to $\Delta_0$-formulae in weak arithmetics.
 
Now we need need to define two important global functions we use in local recursions. The first is the \emph{concatenation function} on linear orderings of disjoint linear orderings.

\begin{definition}[Concatenation function]\label{concat} Let $L = [L_1, \cdots, L_n]$ be a linear ordering of linear orderings $L_1, \cdots, L_n$ such that $\Field(L_i) \cap \Field(L_j) = \varnothing$ when $i \neq j$. Then  $\mathrm{Conc}(L)$ \emph{(}or $L_1 * \cdots * L_n$\emph{)} is defined to be the linear ordering such that
\begin{enumerate}
\item[\emph{(i)}] $\Field(L) = \bigcup\{\Field(L_1), \cdots, \Field(L_n)\}$; and
\item[\emph{(ii)}] If $a, b \in \Field(L)$, then $a <_\mathrm{Conc(L)} b$ if, and only if,
\begin{enumerate}
\item[\emph{(a)}] $a \in L_i$ and $b \in L_j$ and $L_i <_L L_j$; or
\item[\emph{(b)}] $a, b \in L_i$ and $a <_{L_i} b$.
\end{enumerate}
\end{enumerate}
\end{definition}
The second is the \emph{local rank function}, which we need to define the \emph{rank} of a set $S$.
 
\begin{definition}[Local rank function]\label{rankdef}Let $S$ be a given set. First, for each linear ordering $L$ in the epsilon fan of $S$, define
\[
[L] = \{L' \in \epsilon(S) : L \simeq_\mathrm{c} L'\}.
\]
Next, define
\[
\hat S = \{[L] : L \in \epsilon(S)\}
\]
And, finally, define $\Rank(S)$ to be the linear ordering on $\hat S$ induced by the relation: $[L_1] < [L_2] \equiv_{\mathrm{df.}}  L_1 <_\mathrm{c}  L_2$.
\end{definition}

This local rank function has the usual properties of the rank function in infinitary ZF, except that, naively, the rank of a set is defined only up to equivalence modulo $\simeq_\mathrm{o}$, so that sets $S$ and $T$ are `equal in rank' if, and only if, $\mathrm{Rank}(S)\simeq_\mathrm{o}\mathrm{Rank}(T)$.\footnote{We need to define the rank of a set \emph{locally} in this fashion, because, unlike $ZF$, $EA_0$ does not have a natural number sequence available.} In particular, it is easy to show that
\begin{enumerate}
\item[(i)] $\Rank(\varnothing) = [\, ]$;
\item[(ii)] $\Rank(\P(a)) \simeq_\mathrm{o} \Rank(a) +_\mathrm{c} \mathbf{1}$ (where $+_\mathrm{c}$ is the cardinal addition function defined below.)
\item[(iii)] If $a\in S$, then $\Rank(a) +_\mathrm{c} \mathbf{1} \leq_\mathrm{o}  \Rank(S)$.
\end{enumerate}

\section{Arithmetical functions and relations}\label{arithfuncrel}

One way to recapture conventional natural number arithmetic in Euclidean arithmetic is to introduce the notion of a global function being \emph{arithmetical} in the sense that the cardinality of its value depends only on the cardinalities of its arguments.\footnote{Another way is to introduce a finitary version of Dedekind's \emph{simply infinite systems}. This we shall do in Section \ref{sis}.}    

\begin{definition}[Arithmetical global function]\label{arithglofunc} Let $\varphi$ be a global function of one argument. Then $\varphi$ is \emph{arithmetical} if, and only if,
\[
\forall x\forall y[x\simeq_\mathrm{c}  y \supset \varphi(x) \simeq_\mathrm{c} \varphi(y)]
\]
\emph{(}and similarly for global functions of more than one argument\emph{).} 
\end{definition}




Dealing with global arithmetical functions means, in effect, that we are doing \emph{cardinal arithmetic}, except that the `cardinal numbers' here are not objects but infinite equivalence classes of the global relation of cardinal equivalence ($\simeq_\mathrm{c} $).  \bigskip

\noindent \textbf{Examples of Arithmetical Global Functions and Relations}
\begin{enumerate}
\item[(i)]Cardinal successor: $\mathrm{Succ}_\mathrm{c}(a) =_{\mathrm{df.}}  a \cup \{a\}$.

Given a classical natural number $\mathbf{n}$, that belongs to our classical metatheory, we abuse notation and write $\mathbf{n}$ for the constant term $\mathrm{Succ}^{\mathbf{n}}_\mathrm{c}(\varnothing)$, where $\mathrm{Succ}^{\mathbf{n}}_\mathrm{c}(a) =_{\mathrm{df.}}  \underbrace{\mathrm{Succ}(\mathrm{Succ}(\cdots \mathrm{Succ}(a) \cdots))}_{\mathbf{n}}$.  Thus, we use the symbol `$\mathbf{n}$' to denote classical natural numbers in the metatheory and also constant terms of our theory, $EA$, which name sets containing that many elements.  We will ensure that no ambiguity results.
\item[(ii)]Cardinal addition: $ a +_\mathrm{c} b =_{\mathrm{df.}}  (a \times \{\mathbf{0}\}) \cup (b \times \{\mathbf{1}\})$.
\item[(iii)]Cardinal multiplication: $ a \times_\mathrm{c} b =_{\mathrm{df.}}  a \times b$.
\item[(iv)]Cardinal exponentiation: $\exp_\mathrm{c}(a,b) =_{\mathrm{df.}}  a^b\,\, (= \{f \subseteq b \times a : f: b \rightarrow a \})$

Using definition by (external) meta-theoretical recursion, let $\mathbf{2}^a_{\mathbf{0}} = a$ and $\mathbf{2}^a_{\mathbf{k+1}} = \mathbf{2}^{\mathbf{2}^a_{\mathbf{k}}}$.  And we write $\mathbf{2_k}$ for $\mathbf{2^1_k}$.
\item[(v)]Cardinal equality: $a\simeq_\mathrm{c} b$
\item[(vi)]Cardinal inequality: $a\leq_\mathrm{c} b$ 
\end{enumerate}
It is also possible to define\footnote{See \cite{jm}, Section 9.1 for details.} bounded sums  ($\sum_{x<_\mathrm{c} S}\varphi(x)$) and products ($\prod_{x<_\mathrm{c} S}\varphi(x)$), and arithmetical bounded quantifiers of both sorts
\[
(\forall x\leq_\mathrm{c} S)\Psi(x)\mbox{ and }(\exists x\leq_\mathrm{c} S)\Psi(x)
\]
All four of these definitions make use of a unary arithmetical global function, $S\mapsto\mathrm{Card}(S)$, which has the following two key properties:
\begin{enumerate}
\item[(a)] $\mathrm{Card}(S)\simeq_\mathrm{c} S$.
\item[(b)]For all sets $T$, if $T<_\mathrm{c} S$, then $(\exists x\in\mathrm{Card}(S))[x\simeq_\mathrm{c} T]$. 
\end{enumerate}
$\mathrm{Card}(S)$ is called the \emph{local cardinal of} $S$.\footnote{$\mathrm{Card}(S)$ is defined as follows. Let $\hat{S}$ be the set of all equivalence classes of \emph{proper} subsets of $S$ under the $\simeq_\mathrm{c}$.  Linearly order these equivalence classes by the sizes of their members to obtain a linear ordering $L_{\hat{S}}$, and define $\mathrm{Card}(S)$ to be the set of all \emph{proper} initial segments of $L_{\hat{S}}$.
} The \emph{local ordinal}, $\mathrm{Ord}(S)$, of $S$ can then be defined to be the linear ordering of $\mathrm{Card}(S)$ which arranges its members in order of increasing cardinality.

Using the local cardinal function we can add three further arithmetical global functions to our list.
\begin{enumerate}
\item[(vii)] Given a classical natural number, $\mathbf{n}$: $\sqrt[\mathbf{n}]{a} =_{\mathrm{df.}}  \mathrm{min}_\mathrm{c}(\{x \in \mathrm{Card}(a) : a \leq_\mathrm{c} x^{\mathbf{n}}\})$.  Thus, for any set $a$, $\sqrt[\mathbf{n}]{a^{\mathbf{n}}} \simeq_\mathrm{c} a$, but $a \leq_\mathrm{c} (\sqrt[\mathbf{n}]{a})^{\mathbf{n}}$.
\item[(viii)]Base $S$ logarithm: $\mathrm{log}_S(a) =_{\mathrm{df.}}  \mathrm{min}_\mathrm{c}(\{x \in \mathrm{Card}(a) : a \leq_\mathrm{c} S^x\})$.  Thus, for any set $a$, $\mathrm{log}_S(S^a) \simeq_\mathrm{c} a$, but $a \leq_\mathrm{c} S^{\mathrm{log}_S(a)}$.
\item[(ix)] It is a consequence of the Global Function Bounding Lemma (Meta{-}lemma \ref{termbndlem}, proved below) that we cannot define a two place function, $\lambda xy\, \mathbf{2}^x_y$, that satisfies the following recursion equations: $\mathbf{2}^a_{\varnothing} = a$; and $\mathbf{2}^a_{b + \mathbf{1}} = \mathbf{2}^{\mathbf{2}^a_b}$. However, we can define the graph of such a function:  that is, we can define a global relation $\Phi$ of three places such that, for every $a$, $\Phi(a, \varnothing, a)$ and, if $\Phi(a, b, c)$, then $\Phi(a, b+\mathbf{1}, \mathbf{2}^c)$.

By definition, $\Phi(a, b, c)$ if, and only if, there is a linear ordering $L$ of a subset of $\mathrm{Card}(c + \mathbf{1})$ such that
\begin{enumerate}
\item[(a)] $\First(L) \simeq_\mathrm{c} a$;
\item[(b)] For all $x$ in $\Field(L)$ except $\Last(L)$, $\Next_L(x) \simeq_\mathrm{c} \mathbf{2}^x$;
\item[(c)] $\Last(L) \simeq_\mathrm{c} c$; and
\item[(d)] $L \simeq_\mathrm{c} b + \mathbf{1}$.
\end{enumerate}
We write $\mathbf{2}^a_b \simeq_\mathrm{c} c$ for $\Phi(a, b, c)$.  Similarly, we can define the relations $\mathbf{2}^a_b \leq_\mathrm{c}  c$ and $c \leq_\mathrm{c} \mathbf{2}^a_b$.  As above, when $a = \mathbf{1}$, we simply write $\mathbf{2}_b \simeq_\mathrm{c} c$, $\mathbf{2}_b \leq_\mathrm{c}  c$, and $c \leq_\mathrm{c} \mathbf{2}_b$ respectively.  These allow us to define the following arithmetical global function:
\[
\mathrm{suplog}_{\mathbf{2}}(a) =_{\mathrm{df.}} \mathrm{min}_\mathrm{c}(\{x \in \mathrm{Card}(a): a \leq_\mathrm{c} \mathbf{2}_{x}\})
\]
\end{enumerate} 


\section{$EA$ and its relation to fragments of $PA$}\label{EAweakarith}

In the previous section, we described a way of representing arithmetic in Euclidean Arithmetic.  In this section, we consider the relations of interpretability that hold between $EA_1$ and fragments of Peano arithmetic.  We will state the salient theorems, but we will not give detailed proofs.  Our purpose, as we will explain at the end of the section is to prevent a possible misunderstanding of the status of results later in the paper.

In his doctoral thesis \cite{vh}, V. Homolka proved the following theorem:
\begin{theorem}[Homolka]
$EA_1$ and $I\Delta_0 + \mathrm{exp}$ are mutually interpretable.
\end{theorem}
\emph{Proof sketch}:  The interpretation, $\dag$, of $I\Delta_0 + \mathrm{exp}$ in $EA_1$ maps `0' to `$\varnothing$', `$a = b$' to `$a \simeq_\mathrm{c} b$', `$a < b$' to `$a <_\mathrm{c} b$', and the successor, addition, multiplication, and exponentiation function symbols are mapped to the corresponding arithmetic global functions of $EA_1$.

In the other direction, the intepretation, $*$, of $EA_1$ in $I\Delta_0 + \mathrm{exp}$ is the well-known Ackermann interpretation based on the following interpretation of the membership relation:\medskip

$(n \in m)^* \equiv_{\mathrm{df.}.}$ the $n^\mathrm{th}$ bit of $m$ is 1. \hfill $\Box$\medskip

However, it is clear that these interpretations are not inverses of one another.  Given a formula $A$ of $I\Delta_0 + \mathrm{exp}$, $(A^\dag)^*$ is not $A$; and, given a formula $B$ of $EA_1$, $(B^*)^\dag$ is not $B$.  

What's more, neither $\dag$ nor $*$ have inverses.  In the case of $\dag$, this is obvious.  It is not so clear in the case of $*$.  For instance, consider the usual Ackermann interpretation of `ZF with the axiom of infinity negated' in PA:  as proved in \cite{rktlw}, this has an inverse (providing the set theory is suitably formulated). The reason for the failure in our case is this.  Let $WHP$ (Weak Hiearchy Principle) be the following sentence of $EA_1$:
\[
\forall x\exists y\forall z[z \in y \equiv \Rank(z) \leq_\mathrm{c} \Rank(x)]
\]
Thus, $WHP$ asserts the existence of a conventional rank function:  from every set, you can recover the first level of the cumulative hierarchy at which it appears.  Now, it is easy to see that
\[
I\Delta_0 + \mathrm{exp} \vdash (WHP)^*
\]
Under the Ackermann coding, the first level of the cumulative hierarchy at which the set coded $m$ occurs is $2_n - 1$, where $2_{n-2} \leq m < 2_{n-1}$.  Since it follows that $2_n - 1 <2^{2^m}$, a bounded induction establishes the existence of $2_n - 1$, and thus $(WHP)^*$.

However, it is a corollary of Theorems \ref{VNZconstant}, \ref{nounivscale}, and \ref{noACK0pred} below that
\[
EA_1 \not \vdash WHP
\]
Each of those theorems exhibits a formula that $EA_1$ cannot prove; if $EA_1$ were to prove $WHP$, then it would be able to prove those formulae. This is enough to show that $*$ does not have an inverse.

This raises an obvious question:  Does the Ackermann interpretation of $EA_1 + WHP$ in $I\Delta_0 + \mathrm{exp}$ have an inverse?  The answer is yes.  This is Theorem 12.2.2 of \cite{rp}.

These results have some importance for the status of certain of our results.  Given the mutual interpretability of $EA_1$ and $I \Delta_0 + exp$, and given the considerable literature on set theory carried out in $I\Delta_0 + \mathrm{exp}$ via the Ackermann coding, it might be thought that many of our theorems are merely restatements of results already known.   But this is a mistake.  As we saw above, the Ackermann interpretation of $EA_1$ in $I\Delta_0 + \mathrm{exp}$ does not have an inverse:  thus, a set-theoretical result, coded by the Ackermann coding and proved in $I\Delta_0 + \mathrm{exp}$, may not be provable in $EA_1$.  $WHP$ is the obvious example of such a formula, but Theorems \ref{VNZconstant}, \ref{nounivscale}, and \ref{noACK0pred} provide further instances of the phenomenon.  What's more, even when a set-theoretical construction is possible in $EA_1$ and (Ackermann coded) in $I\Delta_0 + \mathrm{exp}$, the construction is often considerably more difficult to effect in the weaker system of $EA_1$, where $WHP$ is absent:  the construction of $\ACK_0$ in Section \ref{lexordsect} is a particularly vivid instance of this; it would be straightforward to define this system in the presence of $WHP$; in its absence, we require a surprising trick.



\section{Iteration systems and natural number systems}\label{sis}

In Section \ref{arithfuncrel}, we showed how to treat natural numbers as equivalence classes of sets under the equivalence relation, $\simeq_\mathrm{c}$, of cardinal equivalence.  In this section, we present the alternative treatment of arithmetic that is our main concern.  We follow the example of Dedekind and attempt to characterise natural numbers directly as the objects generated from a fixed initial object, $a$, by successive iterations of a suitable \emph{successor function}, $\sigma$. In this way we obtain what Dedekind called a \emph{simply infinite system}:
\[
a,\,\sigma(a),\,\sigma(\sigma(a)), \cdots
\]
Since a simply infinite system contains infinitely many terms we cannot define such a system directly. But we can define the class of all its initial segments
\[
[\,],\,\,[a],\,\,[a,\,\sigma(a)],\,\,\,[a,\sigma(a),\sigma(\sigma(a))], \cdots
\]
Such a class of initial segments is what we shall call a \emph{natural number system}.  We make these informal observations precise in the following definitions.

We begin by defining what it means for a linear ordering to be \emph{generated} from an object $a$ by a global function $\sigma$.

\begin{definition}  Let $\sigma$ be a global function of one argument. Then we say that \emph{$L$ is generated from $a$ by $\sigma$} \emph{(}written $\mathrm{Gen}_{\sigma,a}(L)$\emph{)} if $L$ is a linear ordering and either

\noindent\emph{(a)} $L=[\,]$ or

\noindent\emph{(b)} $L\neq[\,]$ and
\begin{enumerate}
\item[\emph{(i)}]$\First(L)=a$ 
\item[\emph{(ii)}]For all $x$ in $\mathrm{Field}(L)$ except $\Last(L)$, $\Next_L(x)=\sigma(x)$.
\end{enumerate}
We call the class of linear orderings generated from $a$ by $\sigma$ the \emph{iteration system $\N_{\sigma, a}$}: Thus
\[
L\mbox{ in }\N_{\sigma, a}\equiv_{\mathrm{df.}}\mathrm{Gen}_{\sigma, a}(L)
\]
\end{definition}
Note that the \emph{definiens} in the definition of  $\mathrm{Gen}_{\sigma,a}(L)$ here is an $EA_0$ formula. (The same obviously applies to the definition of $L\mbox{ in }\N_{\sigma, a}$.)

The empty ordering $[\,]$ and the one-termed ordering $[a]$ are generated from any $a$ by any $\sigma$, and if $\sigma(a)=a$, these are the only linear orderings so generated. In fact, if $\sigma$ generates a linear ordering $L$ from $a$ for which $\sigma(\Last(L))\in \mathrm{Field}(L)$ then the linear orderings generated from $a$ by $\sigma$ are precisely the initial segments of $L$. If no such linear ordering is generated then $\sigma$ generates a natural number system from $a$.

\begin{definition}\label{sisdef} Let $\N_{\sigma, a}$ be an iteration system. We say that \emph{$\N_{\sigma, s}$ is a natural number system} if, for all linear orderings $L$
\[
L\mbox{ in }\N_{\sigma, a} \supset \sigma(\Last(L))\not\in\mathrm{Field}(L)
\]
\end{definition}
Note that in this definition the \emph{definiens} can be expressed by a $\Pi_1$ formula of  $EA_1$.

Thus natural number systems are iteration systems, $\N_{\sigma, a}$, in which the `process' of `generation' never comes to a halt. The \emph{numbers} of a natural number system $\N_{\sigma, a}$ are the linear orderings that compose it as a class.

As we have already explained, in Euclidean arithmetic these linear orderings correspond to the finite initial segments of the infinite sequence of terms 
\[
a,\,\sigma(a),\,\sigma(\sigma(a))\cdots
\]
in conventional infinitary set theory. We therefore define the \emph{terms} of the natural number system $\N$ to be the terms of the linear orderings that are its \emph{numbers}. Thus
\[
a,\,\sigma(a),\,\sigma(\sigma(a))\cdots
\] 
are the \emph{terms} of the natural number system whose \emph{numbers} are
\[
[\,],\,[a],\,[a,\sigma(a)],\,[a,\sigma(a),\sigma(\sigma(a))],\,\cdots
\]
It follows that the numbers of any natural number system are \emph{always} linear orderings; its terms may, but need not, be.
 
The property of \emph{being a number of $\N_{\sigma, a}$} (\emph{an $\N_{\sigma, a}$-number}) is easily expressed in $EA_0$ by a quantifier-free formula
\[
L\mbox{ is an number of }\N_{\sigma,a} \equiv\mathrm{Gen}_{\sigma, a}(L)
\]
whereas, to express the property of \emph{being a term of $\N_{\sigma,a}$} (\emph{an $\N_{\sigma, a}$-term}), we must turn to the first-order theory, $EA_1$, where this property is $\Sigma_1$:
\[
b\mbox{ is a term of }\N_{\sigma, a} \equiv \exists L[\mathrm{Gen}_{\sigma, a}(L) \ \&\  b\in\mathrm{Field}(L)]
\]
Thus the property of being a \emph{number} of $\N_{\sigma, a}$  is  $\Delta_0$ (i.e., expressible as an $EA_0$ formula), whereas that of being a \emph{term} of $\N_{\sigma, a}$ is $\Sigma_1$.

\begin{definition}
Let $\N$ be a natural number system and $\varphi$ a unary global function. We say that $\varphi$ \emph{recovers $\N$-numbers from $\N$-terms} if, for all linear orderings $L$,  
\[
L \mbox{ in } \N \,\supset\, \varphi(\Last(L))= L
\]
\end{definition}

In Section \ref{lexordsect}, we define a natural number system in which it is not possible to define a function in $EA_0$  that recovers its numbers from its terms.  However, in all other systems considered in this paper, it is possible to define such a function.  Whenever it is possible to recover the numbers of a natural number system $\N$ from its terms,  the property of \emph{being a term of $\N$} can be expressed in $EA_0$:
\[
\mbox{$a$ is an $\N$-term} \equiv \mbox{$\varphi(a)$ lies in $\N$}
\]
In these cases, we write $[\cdots, a]_\N$ for $\varphi(a)$ whenever $a$ is an $\N$-term.

In what follows, it will be useful to introduce the following notational convention:  If $\N$ is a natural number system and $\mathbf{k}$ is a classical natural number, by definition,
\[
\mathbf{k}_{\N} =_{\mathrm{df.}} \sigma^{\mathbf{k}}_{\N}(0_{\N})
\]
That is, $\mathbf{k}_{\N}$ is the $\mathbf{k}^{\mathrm{th}}$ term of $\N$, where  $0_\N$ is the initial term of $\N$ and $\sigma_\N$ is its successor function.

Let us consider some examples.\medskip

\noindent \textbf{Examples of Natural Number Systems} 
 
\begin{enumerate}
\item[(i)] The \emph{von Neumann} natural number system, $\VN$:
\begin{enumerate}
\item[] $0_{\VN}=\varnothing$;\ \ \ $\sigma_{\VN}(x)=x\cup\{x\}$
\end{enumerate}
\item[(ii)] The \emph{Zermelo} natural number system, $\Z$:
\begin{enumerate}
\item[]$0_\Z=\varnothing$;\ \ \ $\sigma_\Z(x)=\{x\}$
\end{enumerate}
\item[(iii)] The \emph{cumulative hierarchy} natural number system, $\CH$:
\begin{enumerate}
\item[] $0_\CH=\varnothing$;\ \ \ $\sigma_\CH(x)=\P(x)$
\end{enumerate}

\end{enumerate}

Notice that in all the examples given, we can recover the numbers of the system from its terms:  that is, for each of these three natural number systems, we can define a global function whose value at a term, $b$, of that system is the linear ordering $[0_\N\cdots, b]_\N$ that is the \emph{number} of the system whose last term is $b$. This follows from the fact that in the von Neumann system and the Cumulative Hierarchy system every term preceding $b$ lies in $b$, and in the Zermelo systems they lie in $\mathrm{TC}(b)$.


We can establish mathematical induction in natural number systems, as a derived rule of inference in $EA_0$ and as a provable sentence in $EA_1$.
In fact, we can prove mathematical induction for \emph{arbitrary} iteration systems $\N$, whether or not they are natural number systems.

\begin{theorem}[Mathematical induction in $EA_0$]\label{ea0mathindN} Let $\N$ be an iteration system and $\Phi$ be a global relation. Then from the premises
\begin{enumerate}
\item[\emph{(i)}]$EA_0\vdash \Phi([\, ])\ \&\ \Phi([0_\N])$.
\item[\emph{(ii)}]$EA_0\vdash [0_\N, \cdots, b] \mbox{ in } \N\ \&\ \Phi([0_\N, \cdots, b])\ . \supset .\ \Phi([0_\N, \cdots, b, \sigma_\N(b)])$.
\end{enumerate}
we may infer the conclusion
\[
EA_0\vdash L\mbox{ in }\N \supset \Phi(L)
\]

\end{theorem}
Note that there are no restrictions on the $EA_0$ global relation $\Phi$ (other than that it \emph{be} $EA_0$).
The formulation of mathematical induction in $EA_1$ is equally straightforward.

\begin{theorem}[Mathematical induction in $EA_1$]\label{ea1mathindN} Let $\N$ be an iteration system and $\Phi$ be a quantifier-free global relation.  Then
\medskip

$EA_1\vdash \Phi([\, ])\ \&\ \Phi([0_\N])\ \&\ (\forall x)[x \mbox{ in } \N\ \&\ \Phi(x)\ . \supset .\ \Phi(x * [\sigma_\N(\Last(x))])]$

\hspace{70mm} $\ . \supset .\ (\forall x)[x \mbox{ in } \N \supset \Phi(x)]$

\end{theorem}
Note that this holds for arbitrary iteration systems, not just for natural number systems.

In Section \ref{EAweakarith}, we tried to prevent a possible misunderstanding of the relationship between $EA_1$ and systems of weak arithmetic.  Here, we emphasize another disanalogy.  One might compare our \emph{natural number systems} with \emph{cuts} in models of arithmetic.  After all, they are classes containing an initial element and closed under a successor function.  Thus, when we prove below the existence of natural number systems with various properties, one might be tempted to read these results as restatements of results concerning models of arithmetic that assert the existence of cuts with analogous properties.

This would be a mistake.  Given a model of $EA_1$, the numbers of a particular natural number system in that model will provide a model of some weak system of arithmetic; but they are not \emph{cuts} in some larger model of arithmetic.  In particular, for the same reasons as given in Section \ref{EAweakarith}, it is not the case that every model of $EA_1$ can be transformed into a model of $I\Delta_0 + \mathrm{exp}$ via the Ackermann interpretation.  And, even in those models in which this is possible, the numbers of a given natural number system in that model will not correspond to a cut, since the successor function of the natural number system will not correspond to the successor function in the arithmetic.

\section{Length and  closure for natural number systems}\label{lengthclose}

We are interested in two sorts of question concerning natural number systems: questions concerning comparative length among such systems, and questions concerning their closure under various arithmetical functions. Let us begin by introducing some useful concepts and notation.

\begin{metadefinition} Let $\M$ and $\N$ be natural number systems, $\varphi$ be an arithmetical global function of one argument, and $\eta$ is a global function of one argument. Then \emph{$\eta$ represents $\varphi$ as a map from $\M$ to $\N$} \emph{(}in symbols, $\eta:\M\stackrel{\varphi}{\rightarrow}\N$\emph{)} if
\[
(\forall x)[x \mbox{ in } \M\, .\supset.\, \eta(x) \mbox{ in } \N \ \&\  \varphi(x)\simeq_\mathrm{c}\eta(x)]
\]
Similarly for $\varphi$ and $\eta$ of two or more arguments.\footnote{Recall that for linear orderings $L$, $L\simeq_\mathrm{c}\Field(L)$.}
\end{metadefinition}
Note that this definition is $\Pi_1$. If $\varphi$ is arithmetical we write $\M\stackrel{\varphi}{\rightarrow}\N$ to mean that we can define a global function, $\eta$, such that $\eta:\M\stackrel{\varphi}{\rightarrow}\N$.

As a straightforward consequence of the definition and notational conventions just laid down we know that if $\eta_1:\N_1\stackrel{\varphi}{\rightarrow}\N_2\mbox{ and }\eta_2:\N_2\stackrel{\psi}{\rightarrow}\N_3$, then $\eta_2\circ\eta_1:\N_1\stackrel{\psi\circ\varphi}{\rightarrow}\N_3$.

\begin{metadefinition}\label{sisleq} Let $\M$ and $\N$ be natural number systems.  Then we say that \emph{$\M$ is shorter than, or equal to, $\N$ in length} \emph{(}in symbols,  $\M\preceq \N$\emph{)} if we can explicitly define a global function, $\mu$, such that
\[
\mu : \M \stackrel{\lambda x\, (x)}{\rightarrow} \N
\]
In this case, we say that \emph{$\mu$ is a measure for $\M$ in $\N$}.\footnote{We call this a metadefinition because it cannot be expressed in $EA_1$.}
\end{metadefinition}
If $\M \preceq \N$ and $\N \preceq \M$, we write $\M \cong \N$ and say that \emph{$\M$ and $\N$ are of the same length}.  We write $\M \not \preceq \N$ if it is \emph{not} possible to define a global function $\mu$ such that
\[
\mu : \M \stackrel{\lambda x\, (x)}{\rightarrow} \N
\]

Given a theory $\mathcal{T}$ extending $EA_1$,  we can speak of a natural number system $\M$'s being shorter than or equal to, or not being shorter than or equal to, the natural number system $\N$ \emph{in the theory $\mathcal{T}$}, which we express symbolically by $\M\preceq_\mathcal{T}\N$ and $\M\not\preceq_\mathcal{T}\N$, respectively.  

In Section \ref{complength}, we will prove that, by this definition, two of the three natural number systems, $\VN$, $\CH$, and $\Z$, defined above are \emph{incommensurable} in $EA_1$:  neither is shorter than or equal to the other, that is, $\M\not\preceq_{EA_1}\N$ and $\M\not\preceq_{EA_1}\N$.
\begin{metadefinition}\label{sisclo} Let $\N$ be a natural number system and $\varphi$ an arithmetical global function of one argument. Then we say that \emph{$\N$ is closed under $\varphi$} if we can explicitly define a global function $\eta$ of one argument such that
\[
\eta:\N\stackrel{\varphi}{\rightarrow}\N
\]
\emph{(}and similarly for arithmetical functions $\varphi$ of two or more arguments\emph{)}.
\end{metadefinition}

On the face of it, there is an alternative to Definition \ref{sisleq}, namely 
\begin{metadefinition}\label{sisleqm} Let $\M$ and $\N$ be natural number systems. Then we say that \emph{$\M$ is shorter than, or equal to, $\N$ in length} \emph{(}in symbols,  $\M\preceq \N$\emph{)} if
\[
(\forall x)(\exists y)[x \mbox{ in } \M\, . \supset.\, y \mbox{ in } \N \ \&\  x \simeq_\mathrm{c} y]
\]
\end{metadefinition}
Similarly, the following appears to be an alternative to Definition \ref{sisclo}:
\begin{metadefinition}\label{sisclom}
Let $\N$ be a natural number system and $\varphi$ an arithmetical global function of one argument. Then we say that \emph{$\N$ is closed under $\varphi$} if 
\[
(\forall x)(\exists y)[x \mbox{ in } \N\, .\supset.\, y \mbox{ in } \N \ \&\  y \simeq_\mathrm{c} \varphi(x)]
\]
\end{metadefinition}
In fact, Metadefinition \ref{sisleq} and Definition \ref{sisleqm} are equivalent in $EA_1$, as are Metadefinition \ref{sisclo} and Definition \ref{sisclom}.  This is a consequence of the following metatheorem, which is analogous to a celebrated theorem due to Parikh \cite{rp} concerning bounded arithmetics.
\begin{metatheorem}\label{parikh}
Let $A$ be a $EA_0$-formula whose free variables are from among the distinct free variables $a_1,\cdots,a_n,a_{n+1}=\vec{a},a_{n+1}$.  Then the following two propositions are equivalent:
\begin{enumerate}
\item[\emph{(1)}]There is an $EA_0$-term $t$ with free variables from among $a_1,\cdots,a_n$ such that
\[
EA_0 \vdash A[ t/a_{n+1}] 
\]
\item[\emph{(2)}] $EA_1 \vdash (\forall \vec x)(\exists y)A[\vec x/\vec{a}, y/a_{n+1}]$
\end{enumerate}
\end{metatheorem}
\emph{Proof}: Clearly, (1) implies (2).  The implication from (2) to (1) may be established by modifying Parikh's proof in \cite{rp}. It can also be proved in a manner exactly analogous to the well-known model-theoretic proof via compactness of Parikh's original result---see, for instance \cite{krajicek}, p. 65. \hfill $\Box$\bigskip


\noindent However, since our proofs will proceed by showing the possibility and  impossibility of defining certain global functions in $EA_0$, we take Metadefinitions \ref{sisleq} and \ref{sisclo} as our basic definitions. 


\section{Three Important Meta-theorems}\label{twometathms}

Since many of our key results assert that there do \emph{not} exist global functions of $EA_0$ with certain properties, we need to present three important metatheorems that establish limitations on the possibilities for defining global functions in $EA_0$.

The first of these tells us that, for each global function of $EA_0$, we can establish a certain kind of bound for the action of that function.

\begin{metatheorem}[Global Function Bounding Lemma]\label{termbndlem}
Let $t$ be an $EA_0$ term all of whose free variables are from among the distinct variables $\vec{a}=a_1,\,\cdots,\,a_n$. Then we can find a classical natural number $\mathbf{k}$ such that 
\[
EA_1\vdash \forall \vec{x}[\,t[\vec{x}/\vec{a}]\in \P^{\mathbf{k}}(\mathrm{TC}(\{x_\mathbf{1}, \cdots, x_\mathbf{n}\}))\,]
\]
\end{metatheorem}
\noindent\emph{Proof}: The proof is by meta-theoretic induction on the construction of the term $t$.

First, we note the following facts, which can be shown by meta-theoretic induction to be provable in $EA_0$:
\begin{enumerate}
\item[(A)]If $S$ is a transitive set (i.e., $(\forall x\in S)[x\subseteq S]$), then $\P^{\mathbf{n}}(S)$ is also transitive, for all $\mathbf{n}$.
\item[(B)] If $S$ and $T$ are both transitive, then for all  $\mathbf{m} \leq \mathbf{n}$,
\[
\P^{\mathbf{m}}(S) \subseteq \P^{\mathbf{n}}(S \cup T).
\] 
\end{enumerate}

Now we can proceed with the induction on the construction of $t$.\bigskip

\noindent\textbf{Basis Case ($t$ is a free variable or the constant $\varnothing$):}\medskip
 
\noindent If $t$ is a free variable, $a$, and $a$ is amongst $\vec a$, then
\[
EA_0 \vdash t = a \in \mathrm{TC}(\{a\}) = \P^{\mathbf{0}}(\mathrm{TC}(\{a\})) \subseteq \P^{\mathbf{0}}(\mathrm{TC}(\{\vec a\})) 
\] 

\noindent If $t$ is the constant term, $\varnothing$, then
\[
EA_0 \vdash  t= \varnothing \in \P^\mathbf{1}(\mathrm{TC}(\varnothing)) \subseteq \P^{\mathbf{1}}(\mathrm{TC}(\{\vec a\})) 
\] 
\noindent\textbf{Induction Case ($t$ is a complex term built up from simpler terms):}\medskip

\noindent We shall consider only representative subcases here.\bigskip

\noindent\textbf{Subcase ($t=\{t_1,t_2\}$):} By hypothesis there are $\mathbf{k_1}$ and $\mathbf{k_2}$ such that
\[
EA_0 \vdash t_1 \in \P^{\mathbf{k_1}}(\mathrm{TC}(\{\vec{a}\})) \mbox{\ \ and\ \ } EA_0 \vdash t_2 \in \P^{\mathbf{k_2}}(\mathrm{TC}(\{\vec{a}\})).
\]

\noindent Then (by (B) above)
\[
EA_0\vdash \{t_1, t_2\} \in \P^{\max(\mathbf{k_1},\mathbf{k_2})\mathbf{+1}}(\mathrm{TC}(\{\vec{a}\}))
\]

\noindent\textbf{Subcase ($t=\{x\in t_1: A[x/a]\}$):} By hypothesis, there is $\mathbf{k}$ such that
\[
EA_0 \vdash t_1 \in \P^{\mathbf{k}}(\{\mathrm{TC}(\vec{a})\})
\] 
Hence
\[
EA_0\vdash\{x\in t_1:A[x/a]\}\in\P^{\mathbf{k+1}}(\{\mathrm{TC}(\vec{a}\}))
\]

\noindent \textbf{Subcase ($t=\{t_1[x/a]:x\in t_2\}$):} By hypothesis there are $\mathbf{k_1}$ and $\mathbf{k_2}$ such that 
\begin{enumerate}
\item[(i)]$EA_0 \vdash t_1 \in \P^{\mathbf{k_1}}(\mathrm{TC}(\{a,\vec{a}\}))$,
 where $a$ does not occur in the list  $\vec{a}$.
\item[(ii)]$EA_0 \vdash t_2 \in \P^{\mathbf{k_2}}(\mathrm{TC}(\{\vec{a}\}))$.
\end{enumerate} 
We may conclude that
\[
EA_0 \vdash(\forall x\in t_2)[ t_1[x/a] \in \P^\mathbf{k_1 + k_2}(\mathrm{TC}(\{x,\vec{a}\}))]
\]
From this it follows that
\[
EA_0\vdash\{t_1[x/a]:x\in t_2\} \in \P^\mathbf{k_1 + k_2}(\mathrm{TC}(\{\vec{a}\}))
\]

\noindent\textbf{Subcase $t=\epsilon(s)$:} By hypothesis there is $\mathbf{k}$ such that
\[
EA_0 \vdash s \in \P^\mathbf{k}(\mathrm{TC}(\{\vec{a}\}))
\]
it follows that
\[
EA_0 \vdash \epsilon(s) \in \P^\mathbf{k+3}(\mathrm{TC}(\{\vec{a}\}))
\]

\noindent This completes our induction and the lemma follows. \hfill $\Box$ \medskip

Notice that in both the statement and the proof of the Global Function Bounding Lemma, $\mathbf{k}$ is a classical meta-variable ranging over the classical natural numbers, and $t$ is a classical meta-variable ranging over the infinite set of $EA_0$ terms. The induction is a classical induction proof. In short, we are using ordinary, conventional, infinitary mathematics in our meta-theory.  

An important corollary of the Global Function Bounding Lemma is the following:
\begin{metatheorem}[Rank Bounding Lemma]
Let $t$ be an $EA_0$ term without second order variables all of whose first order variables are from among the distinct variables $\vec{a}=a_\mathbf{1},\,\cdots,\,a_\mathbf{n}$. Then we can find a classical natural number $\mathbf{k}$ such that 
\[
EA_1\vdash \forall\vec{x}\,[\,\Rank(t[\vec{x}/\vec{a}]) <_\mathrm{c}  \mathrm{max}_\mathrm{c}(\{\Rank(x_\mathbf{1}), \cdots, \Rank(x_\mathbf{n})\}) + \mathbf{k}\,]
\]
\end{metatheorem}
\emph{Proof}:  This follows from the Global Function Bounding Lemma along with the following two facts:
\begin{enumerate}
\item[(a)] $EA_0 \vdash \Rank(\mathrm{TC}(\{a_\mathbf{1}, \cdots, a_\mathbf{n}\})) \simeq_\mathrm{c} \mathrm{max}_\mathrm{c}(\{\Rank(a_\mathbf{1}), \cdots, \Rank(a_\mathbf{n})\})$
\item[(b)] $EA_0 \vdash \Rank(\P(a)) \simeq_\mathrm{c} \Rank(a) +_\mathrm{c} \mathbf{1}$ \hfill $\Box$
\end{enumerate}

Our third key metatheorem imposes a restriction on correlations between terms of the natural number systems, $\VN$ and $\Z$. 

\begin{metatheorem}[$\VN$-$\Z$ Splitting Lemma]\label{vnz}
Suppose $\mathbf{n}$ is a classical natural number.  Then\smallskip

\noindent $EA_1\vdash \forall x\forall y[\mathrm{Term}_\VN(x)\ \&\  \mathrm{Term}_\Z(y)\ . \supset .\ $

\hspace{40mm}$\P^{\mathbf{n}}(\mathrm{TC}(\{x\}))\cap \P^{\mathbf{n}}(\mathrm{TC}(\{y\}))\subseteq \P^{\mathbf{n+2}}(\varnothing)]$\footnote{Recall that in both $\VN$ and $Z$ numbers can be recovered from terms, so that $\mathrm{Term}_\VN$ and $\mathrm{Term}_\Z$ are $EA_0$-formulas.}
\end{metatheorem}

\noindent\emph{Proof}: We proceed by induction on $\mathbf{n}$.\medskip

\noindent\textbf{Basis Case ($\mathbf{n}=\mathbf{0}$)}: Suppose that $v$ is a term of $\VN$, $z$ is a term of $\Z$, and $a \in \mathrm{TC}(\{v\}) \cap \mathrm{TC}(\{z\})$.  Then, since every element of $\mathrm{TC}(\{v\})$ is a term of $\VN$ and every term of $\mathrm{TC}(\{z\})$ is an element of $\Z$, and since the only terms of $\VN$ that are also terms of $\Z$ are $\varnothing$ and $\{\varnothing\}$, it follows that $a=\varnothing$ or $a=\{\varnothing\}$, and thus, $a\in\P^\mathbf{2}(\varnothing)$, as required.\medskip

\noindent\textbf{Induction Case ($\mathbf{n}>\mathbf{0}$)}: Suppose as hypothesis of induction, that if $v$ is a term of $\VN$ and $z$ is a term of $\Z$, then
\[
\P^{\mathbf{n-1}}(\mathrm{TC}(\{v\})) \cap \P^{\mathbf{n-1}}(\mathrm{TC}(\{z\}))\subseteq\P^{\mathbf{n+1}}(\varnothing)
\]
we must prove that
\[
\P^{\mathbf{n}}(\mathrm{TC}(\{v\})) \cap \P^{\mathbf{n}}(\mathrm{TC}(\{z\}))\subseteq\P^{\mathbf{n+2}}(\varnothing)
\]
To that end, let $a\in \P^{\mathbf{n}}(\mathrm{TC}(\{v\})) \cap \P^{\mathbf{n}}(\mathrm{TC}(\{z\}))$ then, every member of $a$ is a member of $\P^{\mathbf{n-1}}(\mathrm{TC}(\{v\})) \cap \P^{\mathbf{n-1}}(\mathrm{TC}(\{z\}))$ and therefore, by hypothesis of induction, a member of $\P^{\mathbf{n+1}}(\varnothing)$. Hence $a\subseteq\P^{\mathbf{n+1}}(\varnothing)$, and therefore $a\in\P^{\mathbf{n+2}}(\varnothing)$, for every $a\in \P^{\mathbf{n}}(\mathrm{TC}(\{v\})) \cap \P^{\mathbf{n}}(\mathrm{TC}(\{z\}))$. It follows that
\[
\P^{\mathbf{n}}(\mathrm{TC}(\{v\})) \cap \P^{\mathbf{n}}(\mathrm{TC}(\{z\}))\subseteq\P^{\mathbf{n+2}}(\varnothing)
\]
 as required. This establishes the induction case and the theorem follows.\hfill$\Box$ \medskip

As we shall see, the $\VN$-$\Z$ Splitting Lemma has some important consequences for the general theory of natural number systems.\bigskip

\section{Comparing Lengths among  $\VN$, $\Z$, and $\CH$}\label{complength}

Let us now turn our attention to the three natural number systems $\VN$ (the \emph{von Neumann} system), $\Z$ (the \emph{Zermelo} system), and $\CH$ (the \emph{Cumulative Hierarchy} system) introduced in Section \ref{sis}. First notice that $\CH\preceq\VN$ and $\CH\preceq\Z$.  For given any linear ordering, $L$, that is a number of $\CH$, the set of terms of $\VN$ that lie in $\P(\Last(L))$ are the terms of a linear ordering that is a number of $\VN$ of the same length as $L$. This shows that $\CH\preceq\VN$; the argument that $\CH\preceq\Z$ is essentially the same.

The remaining question concerning the length relations between $\VN$ and $\Z$ was answered by S. Popham in \cite{sp}: they are incommmensurable in $EA$. Popham's argument is model theoretic; ours is related to his, since both appeal to the $\VN$-$\Z$ Splitting Lemma, but ours is more syntactic in character.

\begin{metatheorem}[Popham]\label{popthm}$\VN$ and $\Z$ are incommensurable in $EA_1$:  that is,
\[
\VN\not\preceq_{EA_1} \Z\mbox{ and }\Z\not\preceq_{EA_1} \VN
\]
\end{metatheorem}
\emph{Proof}: We derive this as a straightforward consequence of the following lemma. \hfill $\Box$\medskip 

\begin{metatheorem}\label{VNZconstant}Let $\varphi$ be of the form $\lambda x\, t[x / a]$ where $t$ contains no free variables other than $a$.  Then
\begin{enumerate}
\item[\emph{(i)}]If $ \varphi :\VN\to\Z$, then $\varphi$ is bounded above by a constant function: i.e., there exists a classical natural number $\mathbf{k}$ 
\[
EA_1\vdash \forall X[X \mbox{ in } \VN \supset \varphi(X) <_\mathrm{c}  \mathbf{k}]
\]
\item[\emph{(ii)}]If $\varphi :\Z\to\VN$, then then $\varphi$ is bounded above by a constant function: i.e., there exists a classical natural number $\mathbf{k}$ such that 
\[
EA_1 \vdash \forall X[X \mbox{ in } \Z \supset \varphi(X)  <_\mathrm{c}  \mathbf{k}]
\]
\end{enumerate}
\end{metatheorem}
\noindent\emph{Proof}:  We shall just prove (i) since the proof of (ii) is essentially the same.

Suppose that 
\[
EA_1 \vdash \forall X[X \mbox{ in } \VN \supset \varphi(X) \mbox{ in } \Z]
\]
and suppose $EA_0\vdash \varphi([\,]) = [\cdots, (\mathbf{k_1})_{\Z}]_{\Z}$.  (Recall:  If $\N$ is a natural number system and $\mathbf{n}$ is a (classical) natural number, then by definition $\mathbf{n}_{\N} =_{\mathrm{df.}} \sigma^{\mathbf{n}}_{\N}(0_{\N})$.)

Then define a unary global function $\varphi' = \lambda x\, t'[x / a]$ such that $t'$ contains no free variables other than $a$ and
\begin{enumerate}
\item[(a)] $\forall x[x\mbox{ is not a term of }\VN\supset \varphi'(x) =\varnothing]$.\medskip
\item[(b)] $\forall x[x\mbox{ is a term of }\VN\ \&\ \varphi([\cdots, x]_\VN) \neq [\,] .\supset. $

\hspace{60mm}$\varphi'(x)= \Last(\varphi([ \cdots, x]_\VN))]$.\medskip
\end{enumerate}

\noindent Thus, $\varphi'$ takes terms of $\VN$ to terms of $\Z$ such that, if $L \neq [\, ]$ is a number of $\VN$ then $\varphi(L) = [\, ]$ or $\varphi(L)$ is the number of $\Z$ whose last element is $\varphi'(\Last(L))$:  that is, $\varphi(L) = [\cdots, \varphi'(\Last(L))]_\Z$.

Thus, to prove this lemma, it suffices to show that there is a classical natural number $\mathbf{k}$, such that 
\[
EA_1 \vdash \forall x[\mbox{$x$ is a term of $\VN$} \supset \varphi'(x) \in \{\mathbf{0}_{\Z}, \mathbf{1}_{\Z}, \cdots, \mathbf{k}_{\Z}\}]
\]
For, if this is the case, and $L$ is a number of $\VN$,
\begin{eqnarray*}
\varphi(L) & = & [\, ], [\cdots, \varphi'(\Last(L))]_{\Z}, \mbox{ or } \varphi([\, ]) \\
&= &[\, ],\,\,[ \mathbf{0}_{\Z}],\,\, [\mathbf{0}_{\Z}, \mathbf{1}_{\Z}],\,\, \cdots, [ \mathbf{0}_{\Z}, \cdots, \mathbf{k}_{\Z}],\,\, \mbox{ or }\,\,[ \mathbf{0}_{\Z}, \cdots, (\mathbf{k_1})_{\Z}] \\
& <_\mathrm{c}  & \mathrm{max}_\mathrm{c}(\{\mathbf{k+2}, \mathbf{k_1+2}\})
\end{eqnarray*}
as required.

Now, by the Global Function Bounding Lemma, there is $\mathbf{k_0}$ such that 
\[
EA_1 \vdash \forall x[\varphi'(x) \in \P^\mathbf{k_0}(\mathrm{TC}(\{x\}))].
\]
Also
\[
EA_1 \vdash \forall x[x \in \P^{\mathbf{k_0}}(\mathrm{TC}(\{x\}))]
\]
So, in particular
\[
EA_1 \vdash \forall x[\varphi'(x)\in \P^{\mathbf{k_0}}(\mathrm{TC}(\{\varphi'(x)\}))].
\]
Thus, putting these together we get
\[
EA_1 \vdash \forall x[\varphi'(x) \in \P^\mathbf{k_0}(\mathrm{TC}(\{x\})) \cap \P^{\mathbf{k_0}}(\mathrm{TC}(\{\varphi'(x)\}))].
\]
So, arguing in $EA_1$ using the free variable $a$, if $a$ is a term of $\VN$, and hence $\varphi'(a)$ is a term of $\Z$, then by the $\VN$-$\Z$ Splitting Lemma,
\[
\varphi'(a) \in \P^{\mathbf{k_0 + 2}}(\varnothing)
\]
So, by universal generalization,
\[
EA_1 \vdash \forall x[x \mbox{ is a term of } \VN \supset \varphi'(x) \in \P^{\mathbf{k_0 + 2}}(\varnothing)]
\]
and hence
\[
EA_1 \vdash \forall x[x \mbox{ is a term of } \VN \supset \varphi'(x) \simeq_\mathrm{c} \mathbf{0}_{\Z}, \mathbf{1}_{\Z}, \cdots, \mbox{ or } (\mathbf{k_0 + 1})_{\Z}]
\]
as required. \hfill $\Box$\bigskip

It should be clear that this `failure' result does not hold of the theory $EA_1 + WHP$ considered in Section \ref{EAweakarith}.  Indeed, in the presence of the rank function given by $WHP$, $\VN$ and $\Z$ have the same length:  that is, $\VN \cong_{EA_1 + WHP} \Z$.  The same is true of Metatheorem \ref{nounivscale} below.  In this, the $\VN$-$\Z$ Splitting Lemma is used to show that no natural number system is `long enough' to constitute a \emph{universal scale}.  In $EA_1 + WHP$, this is false; the systems $\ACK$ and $\ACK_0$ of Section \ref{lexordsect} both provide a universal scale in the presence of $WHP$.

\begin{metatheorem}\label{nounivscale} 
Let  $\N$ be a natural number system and let $\varphi$ be of the form $\lambda x\, t[x / a]$ where $t$ contains no free variables other than $a$. Then  $\varphi$ does not provide a measure for the class of all sets relative to $\N$,  that is,
\[
EA_1\not \vdash\forall x[\varphi(x) \mbox{ in } \N \ \&\  \varphi(x) \simeq_\mathrm{c}x\,]
\]
\end{metatheorem}
\noindent\emph{Proof}: Given $\N$ and $\varphi$ as in the premise, suppose, by way of contradiction, that
\[
EA_1\vdash \forall x[\varphi(x) \mbox{ in } \N \ \&\  \varphi(x) \simeq_\mathrm{c}x\,]
\]
Then, by the Global Function Bounding Lemma (Lemma \ref{termbndlem}), there is \emph{fixed} classical natural number $\mathbf{k_0}$ such that
\[
EA_1\vdash \forall x[\varphi(x)\in \P^{\mathbf{k_0}}(\mathrm{TC}(\{x\}))]
\]
Now, for all classical natural numbers $\mathbf{k}$,
\begin{eqnarray*}
\mathbf{k}_{\VN} & = & \{\mathbf{0}_{\VN}, \cdots, (\mathbf{k-1})_{\VN}\} \\
&\simeq_\mathrm{c}& [\cdots, (\mathbf{k-1})_{\N}]_{\N}\\
&\simeq_\mathrm{c}& \mathrm{TC}(\mathbf{k}_{\Z})\\
& =& \{\mathbf{0}_{\Z}, \cdots, (\mathbf{k-1})_{\Z}\}
\end{eqnarray*} 
By assumption, $\varphi$ provides a scale relative to $\N$, so
\[
\varphi(\mathbf{k}_{\VN}) = [\cdots, (\mathbf{k-1})_{\N}]_{\N} = \varphi(\mathrm{TC}(\mathbf{k}_{\Z}))
\]
However,
\begin{enumerate}
\item[]$EA_0 \vdash \varphi(\mathbf{k}_{\VN})\in\P^\mathbf{k_0}(\mathrm{TC}(\{\mathbf{k}_{\VN}\})) \subseteq \P^\mathbf{k_0+1}(\mathrm{TC}(\{\mathbf{k}_{\VN}\}))$
\end{enumerate}
\noindent and
\begin{enumerate}
\item[]$EA_0 \vdash \varphi(\mathrm{TC}(\mathbf{k}_{\Z})) \in\P^\mathbf{k_0}(\mathrm{TC}(\{\mathrm{TC}(\mathbf{k}_{\Z})\}))) \subseteq \P^\mathbf{k_0 + 1}(\mathrm{TC}(\{\mathbf{k}_{\Z}\}))$
\end{enumerate}
Thus,
\begin{enumerate}
\item[]$EA_0 \vdash \varphi(\mathbf{k}_{\VN}) = \varphi(\mathrm{TC}(\mathbf{k}_{\Z})) \in \P^\mathbf{k_0+1}(\mathrm{TC}(\{\mathbf{k}_{\VN}\}))\cap\P^\mathbf{k_0+1}(\mathrm{TC}(\{\mathbf{k}_{\Z}\}))$
\end{enumerate}
\noindent and, therefore, by the $\VN$-$\Z$ Splitting Lemma (Lemma \ref{vnz}), for all $\mathbf{k}$
\begin{enumerate}
\item[] $EA_0 \vdash [\cdots, (\mathbf{k-1})_{\N}]_{\N} \simeq_\mathrm{c} \varphi(\mathbf{k}_{\VN}) \simeq_\mathrm{c}  \varphi(\mathrm{TC}(\mathbf{k}_{\Z}))\in \P^{\mathbf{k_0 + 2}}(\varnothing)$
\end{enumerate}
\noindent But this is obviously impossible, since $\mathbf{k_0}$ is a particular classical natural number, but $\mathbf{k}$ can take on any value whatsoever. The theorem now follows.\hfill$\Box$\bigskip

We have now shown that $\CH\preceq_{EA_1}\VN$, $\CH\preceq_{EA_1}\Z$, $\Z\not\preceq_{EA_1}\VN$, and $\VN\not\preceq_{EA_1}\Z$. $\VN \not\preceq_{EA_1}\CH$ and $\Z \not \preceq_{EA_1} \CH$ follow from the relationships already proved since, if $\VN \preceq_{EA_1} \CH$, then $\VN \preceq_{EA_1}\Z$, which gives a contradiction, and similarly  for $\Z \not \preceq_{EA_1} \CH$. Thus we have established the relationships in respect of length among these three natural number systems in $EA_1$.


\section{The Closure Properties of $\VN$, $\Z$, and $\CH$}\label{clpropvnzch}

Given a natural number system, $\N$, we are interested in the arithmetical global functions, $\varphi$, for which $\N \stackrel{\varphi}{\rightarrow} \N$, where $\varphi$ is a function of one or more arguments.  In fact, since we will primarily be interested in  those arithmetical global functions that are increasing in each argument, we can confine our attention mainly to unary arithmetical global functions when discussing closure properties, as we shall show below.

Given any arithmetical global function $\varphi$ which is increasing in each of its arguments, we can define its \emph{unary reduct}, $\varphi_\mathrm{un}$, to be the unary global function characterized by stipulating that
\[
\varphi_\mathrm{un}(a)=_{\mathrm{df.}}\varphi(a,\cdots,a)
\]
\noindent We can easily show that, if $\varphi$ is arithmetical and increasing in each argument, then
\[
\N \stackrel{\lambda \vec{x}\, \varphi(\vec{x})}{\rightarrow} \N\mbox{ if, and only if }\N \stackrel{\lambda x\, \varphi_\mathrm{un}(x)}{\rightarrow} \N
\]

\noindent Thus, when we are investigating questions of the closure of natural number systems under  arithmetical global functions, $\varphi$, which are increasing in each of their arguments, we can confine our attention to their unary counterparts, $\varphi_\mathrm{un}$.

\begin{proposition}  Given a unary arithmetical global function, $\varphi$,
\[
\N \stackrel{\varphi}{\rightarrow} \N\mbox{ if, and only if, for all natural numbers }\mathbf{n},\,\N \stackrel{\lambda x\, \varphi^{\mathbf{n}}(x)}{\rightarrow} \N
\]
\end{proposition}
\noindent This can be proved by a straightforward (meta-)induction on the classical natural number $\mathbf{n}$.

\begin{definition} Let  $\N$ is a natural number system and  $\varphi$ and $\psi$ be unary arithmetical global functions.
\begin{enumerate}
\item[\emph{(i)}]\emph{$\psi$ eventually dominates $\varphi$} \emph{(}in symbols, $\varphi \ll \psi$\emph{)} if there is a classical natural number, $\mathbf{M}$ such that 
\[
EA_1\vdash\forall x[\mathbf{M} \leq_\mathrm{c}  x \supset \psi(x) <_\mathrm{c}  \varphi(x)]
\]
\item[\emph{(ii)}]  If $\varphi,\psi:\N\to\N$, then \emph{$\varphi$ eventually dominates $\psi$ on $\N$} \emph{(}in symbols, $\psi \ll_{\N}\varphi$\emph{)} if there is a classical natural number, $\mathbf{M}$ such that
\[
EA_1\vdash \forall x[x\mbox{ in } \N \ \&\  \mathbf{M} \leq_\mathrm{c}  x\, .\supset.\,  \psi(x) <_\mathrm{c}  \varphi(x)]
\]
\end{enumerate}
\end{definition}
\begin{proposition} If $\varphi \ll_{\N}\psi$ and $\N \stackrel{\psi}{\rightarrow} \N$, then $\N \stackrel{\varphi}{\rightarrow} \N$.
\end{proposition}

\begin{definition} Let $\N$ be a natural number system and let  $\varphi$ be  a unary arithmetical global function increasing in its argument. Then \emph{$\varphi$ is maximally powerful in $\N$} if 
\begin{enumerate}
\item[\emph{(i)}] $\N \stackrel{\varphi}{\rightarrow} \N$, and 
\item[\emph{(ii)}] For any unary arithmetical global function, $\psi$, such that $\N \stackrel{\psi}{\rightarrow} \N$, there is a number $\mathbf{n}$ such that $\lambda x\, \psi(x) \ll_{\N} \lambda x\,\varphi^{\mathbf{n}}(x)$.
\end{enumerate}
\end{definition}
Thus, given a natural number system, $\N$, we are interested in which arithmetical global functions are maximally powerful in $\N$.  Once we know this, we know under exactly which functions $\N$ is closed:  it is closed under those functions that, on $\N$, are eventually dominated by some finite iterate of a function that is maximally powerful in $\N$.

\begin{theorem}
In $EA_1$, $\lambda x\,(x +_\mathrm{c}\mathbf{1})$ is maximally powerful in $\CH$, $\VN$, and $\Z$.
\end{theorem}
\emph{Proof}:  In each of $\CH$, $\VN$, and $\Z$, the successor function increases rank by one.  By the Rank Bounding Lemma, a global function can only increase rank by a fixed metatheoretical natural number.  From these two facts, the theorem follows.  \hfill $\Box$\medskip

It follows that none of the natural number systems $\CH$, $\VN$, or $\Z$ is closed under addition in $EA_1$.

\section{Hierarchies based on numerals}\label{s-ary}

The closure properties of the natural number systems we have encountered so far have been rather unimpressive.  In each of our three principal examples, the natural number system is closed under cardinal successor and each of its finite iterates but not under anything more powerful.

In this section, we present two ways of defining a hierarchy of natural number systems from a given system, and we will prove that if the initial system has certain properties -- properties that $\CH$, $\VN$, and $\Z$ each possess -- then successive systems in the hierarchy are longer than their predecessors and are closed under more and more powerful arithmetical global functions. 

In the case of both hierarchies, the definition is based on the notion of a \emph{numeral}, which we define first.  For us, a numeral is an ordered pair consisting of a \emph{numeration base} and a \emph{sequence of digits}.
\begin{definition} We say that a set $S$ is a \emph{numeration base} if $S \geq_\mathrm{c} \mathbf{2}$ and $S$ contains exactly one member of each cardinality strictly smaller than its own.
\end{definition}
Thus each von Neumann natural number $ > \mathbf{2}$ is a numeration base; so too is the local cardinal, $\mathrm{Card}(T)$, of any set, $T$, containing two or more elements; also, as will become important in our definition of fixed length numeral hierarchies, given any linear ordering $L$, $\mathrm{InSeg}(L)$ is a numeration base.\footnote{Given a linear ordering $L$, $\mathrm{InSeg}(L)$ is the set of all proper initial segments of $L$. Note that $\mathrm{InSeg}(L)\simeq_\mathrm{c}L$.}

The \emph{digits} (i.e., members) of a numeration base, $S$, are ordered by increasing cardinality so that $0_S$ is always the empty set, $1_S$ has exactly one element, and so on.

As usual, we define a \emph{sequence} to be an ordered pair $s = (L, f)$, where $L (= [0, ..., k])$ is a linear ordering and $f$ is a local function whose domain is $\Field(L)$.  In these circumstances we also say that $s=\<s_0,,\cdots,s_k\>_{[0,\cdots,k]}$ where $s_i$ is $s`i$.
With this terminology in hand, we can give our definition of a numeral:
\begin{definition}
A \emph{numeral} is an ordered pair, $(S, s)$, where $S  = \{0_S, ..., M_S\}$ is a numeration base and $s  = \<s_0, ..., s_k\>_{[0, ..., k]}$ is a sequence such that $s_i \in S$ for $i = 0, ..., k$:  that is, $s$ is a sequence of the digits of $S$, as promised above.\medskip

\noindent We call $[0, ..., k]$ the \emph{length} of the numeral.  A \emph{length $L$ numeral} is a numeral with length $L$.  An \emph{$S$-ary numeral} is a numeral with numeration base $S$. \medskip

\noindent Given a numeral $(S, s)$ and a set $A$, we say that \emph{$(S, s)$ codes the cardinality of $A$} if
\[
A \simeq_\mathrm{c} s_0 \cdot S^{[\, ]} + \cdots + s_k \cdot S^{[0, ..., k]}
\]

\end{definition}
Thus, the order of the digits in our numerals is reversed relative to the usual way of writing them:  the first digit corresponds to units, and subsequent digits to increasing powers of $S$.  For instance, if $S  = \{\mathbf{0}, ..., \mathbf{9}\}$, the cardinality of a set containing two hundred thirty-nine elements would be coded by any  $S$-ary numeral $\<\mathbf{9}, \mathbf{3}, \mathbf{2}\>_L$, where $L$ ia linear ordering of length $\mathbf{3}$. 

Note also that nothing in the definition of numeral just given prevents us from allowing a numeral to end in a tail of zeros.  We will return to this point and its consequences in Section \ref{Saryexthier}.

In each of the two hierarchies of natural number systems that we define in this section, the \emph{terms} of the systems in that hierarchy will be numerals. 

We are interested in ways of generating an unending sequence of numerals.  What's more, we are interested in generating this unending sequence by a successor function that takes a numeral that codes one cardinality to a numeral that codes the next cardinality.  Clearly, there are very many ways to do this.  We will exploit the following two:
\begin{enumerate}
\item[(1)] Fix the numeration base and allow the length of the numeral to increase.

This is the idea behind the definition of the \emph{fixed base} numeral hierarchies presented in Section \ref{Saryexthier}.
\item[(2)] Fix the length of the numeral and allow the size of the numeration base to increase.

This is the idea behind the definition of the \emph{fixed length} numeral hierarchies presented in Section \ref{nthrootext}.
\end{enumerate}

\subsection{Fixed base numeral hierarchies}\label{Saryexthier}

In this section, we exploit the idea in (1) above to define a hierarchy of natural number systems based on a system $\N$ that is already given.  We fix the size of the numeration base $S$ and we allow $L$, the length of the numeral, to increase by letting it range over the numbers of $\N$.  The rest of the section makes this rough recipe precise.

Given a natural number system, $\N$, and a numeration base, $S$, we will define a second natural number system, $\N[S]$, whose \emph{terms} are $S$-ary numerals with lengths in $\N$.\footnote{Recall that an $S$-ary numeral is a numeral with numeration base $S$; and recall that the length of a numeral is the linear ordering that underlies the sequence of digits.} 

Thus, our first job is to define the global function $\mathrm{Num}$, which takes a linear ordering $L$ and a numeration base $S$ to the ordering of all $S$-ary numerals with length an initial segment of $L$, ordered by the cardinality they code.  However, there is a problem.  Given a set $A$, a numeration base $S$, and a linear ordering $L$, there may be two or more $S$-ary numerals with length an initial segment of $L$ that code the cardinality of $A$.  For instance, $\mathbf{3}$ is coded by both of the following two $\{\mathbf{0}, \mathbf{1}\}$-ary numerals with length an initial segment of $[\mathbf{0}, \mathbf{1}, \mathbf{2}, \mathbf{3}]$:
\[
\left ( \{\mathbf{0}, \mathbf{1}\}, \<\mathbf{1}, \mathbf{1}, \mathbf{0}\>_{[\mathbf{0}, \mathbf{1}, \mathbf{2}]}\right ) \mbox{\ \ and\ \ } \left ( \{\mathbf{0}, \mathbf{1}\}, \<\mathbf{1}, \mathbf{1}, \mathbf{0}, \mathbf{0}\>_{[\mathbf{0}, \mathbf{1}, \mathbf{2}, \mathbf{3}]}\right )
\]
Thus, given a numeral $\left (S = \{0_S, ..., M_S\}, \<s_0, ..., s_k\>_{[0, ..., k]} \right )$ we say that it is \emph{proper} if $s_k \neq 0_S$.  Given a set $A$, a numeration base $S$, and a linear ordering $L$, there is at most one proper $S$-ary numeral with length an initial segment of $L$ that codes the cardinality of $A$.  Thus
\begin{definition}
$\mathrm{Num}_S(L)$ is the ordering of the proper $S$-ary numerals with length an initial segment of $L$, ordered by the cardinalities they code.
\end{definition}
Since, for each cardinality, there is at most one $S$-ary numeral with length an initial segment of $L$ that codes that cardinality, this ordering is total.

With this definition in hand, we can define $\N[S]$:
\begin{definition}
$\N[S]$ is generated from the initial element $\left (S, \<0_S\>_{[0_\N]}\right )$ by the following successor function.
\begin{quote}
If $\left (S, s = \<s_0, ..., s_k\>_{[0, ..., k]}\right )$ is an $S$-ary numeral and $[0, ..., k]$ is a number of $\N$, then let\medskip

$\sigma_{\N[S]}\left (\left (S, \<s_0, ..., s_k\>_{[0, ..., k]}\right )\right ) =_\mathrm{df.}$

\hspace{30mm} $\Next_{\mathrm{Num}_S ([0, ..., k, \sigma_\N(k)])}\left (\left (S, \<s_0, ..., s_k\>_{[0, ..., k]}\right )\right )$
\end{quote}
\end{definition}
Clearly, $\N[S]$ is a natural number system.  Further, the \emph{numbers} of $\N[S]$ can be recovered from its \emph{terms} whether or not the same is true of the underlying natural number system, $\N$:  if $\left (S, s = \<s_0, ..., s_k\>_{[0, ..., k]}\right )$ is an $S$-ary numeral and $[0, ..., k]$ is in $\N$, then $[..., (S, s)]_{\N[S]}$ is the initial segment of $\mathrm{Num}_S([0, ..., k])$ that ends with $s$.

We have described a way of producing, given a natural number system, $\N$, and a numeration base, $S$, another natural number system, $\N[S]$.  Metatheoretically,  this process may be iterated indefinitely so that, given  $\N$ and $S$, we can define $\N[S]$, and then $(\N[S])[S]$, and so on.  Let us lay down the following notation by meta-recursion on the natural numbers:
\[
\N^{\mathbf{0}}[S] = \N\,\,\,\,\,\,\,\,\,\,\,\,\,\,\,\N^{\mathbf{k+1}}[S] = (\N^{\mathbf{k}}[S])[S]
\]
We shall call the resulting hierarchy the \emph{hierarchy of $S$-ary numeral extensions of $\N$}.

The question now immediately arises: what are the closure properties of natural number systems that occur in a hierarchy of this sort?  The following lemma and its corollaries are crucial to answering this question.

\begin{lemma}\label{NNS} $\N \stackrel{\lambda x\, S^x}{\rightarrow} \N[S]$ and $\N[S] \stackrel{\lambda x\, \mathrm{log}_S(x)}{\rightarrow} \N$.
\end{lemma}
\noindent \emph{Proof}:  If $S$ is a numeration base and $L$ is a linear ordering, then $\mathrm{Num}_S(L)$ contains exactly one proper $S$-ary numeral for each cardinality $<_\mathrm{c} S^L$.  Thus, $\mathrm{Num}_S(L) \simeq_\mathrm{c} S^L$.  Thus, given a number $L$ in $\N$, let $\eta_1(L) = \mathrm{Num}_S(L)$.  Then $\eta_1 : \N \stackrel{\lambda x\, S^x}{\rightarrow} \N[S]$.\medskip

\noindent By the same reasoning, if $[\cdots, (S, s)]_{\N[S]}$ is in $\N[S]$ and $(S, s)$ has length $[0, ..., k]$, then
\begin{eqnarray*}
S^{[0, ..., k-1]} & \simeq_\mathrm{c} & \left [\cdots, \left (S, \<M_S, ..., M_S\>_{[0, ..., k-1]}\right) \right ]_{\N[S]} \\
& \subsetneqq_* & \left [\cdots, \left (S, s\right) \right ]_{\N[S]} \\
& \subseteq_* & \left [\cdots, \left (S, \<M_S, ..., M_S\>_{[0, ..., k]}\right) \right ]_{\N[S]} \\
& \simeq_\mathrm{c} & S^{[0, ..., k]}
\end{eqnarray*}
Thus,
\[
\mathrm{log}_S\left ( [\cdots, (S, s)]_{\N[S]} \right ) \simeq_\mathrm{c} [0, ..., k]
\]
Thus, let $\eta_1\left ( [\cdots, (S, s)]_{\N[S]} \right ) = [0, ..., k]$.  Then $\N[S] \stackrel{\lambda x\, \mathrm{log}_S(x)}{\rightarrow} \N$. \hfill $\Box$ \medskip

\begin{corollary}\label{cor1NNS}
If $\N \stackrel{\lambda x\, \varphi(x)}{\rightarrow} \N$, then $\N[S] \stackrel{S^{\lambda x\varphi( \mathrm{log}_S(x))}}{\rightarrow} \N[S]$.
\end{corollary}
\noindent\emph{Proof}: By Lemma \ref{NNS}, $\N[S] \stackrel{\lambda x\, \mathrm{log}_S(x)}{\rightarrow} \N \stackrel{\lambda x\, \varphi(x)}{\rightarrow} \N \stackrel{\lambda x\, S^x}{\rightarrow} \N[S]$. \hfill $\Box$

\begin{corollary}\label{cor2NNS}
If $\N[S] \stackrel{\lambda x\, \varphi(x)}{\rightarrow} \N[S]$, then $\N \stackrel{\lambda x\, \mathrm{log}_S(\varphi(S^x))}{\rightarrow} \N$.
\end{corollary}
\noindent\emph{Proof}: By Lemma \ref{NNS}, $\N \stackrel{\lambda x\, S^x}{\rightarrow} \N[S] \stackrel{\lambda x\, \varphi(x)}{\rightarrow} \N[S] \stackrel{\lambda x\, \mathrm{log}_S(x)}{\rightarrow} \N$.
\hfill $\Box$\medskip

\noindent Now we define, by  meta-theoretic recursion,  a hierarchy of arithmetical global functions based on a given arithmetical global function, $\varphi$, as follows:
\[
\varphi^S_{\mathbf{0}}(a) =_{\mathrm{df.}}  \varphi(a) \hspace{10mm} \varphi^S_{\mathbf{k+1}}(a) =_{\mathrm{df.}}  S^{\varphi^S_{\mathbf{k}}(\mathrm{log}_S(a))}
\]
Using this notation, we can prove the following 

\begin{theorem}\label{maxpowNS}
Let $\N$ be a natural number system, $S$ be a numeration base, and $\varphi$ be a unary arithmetical global function that is maximally powerful in $\N$. Then for all $\mathbf{n}$, $\varphi^S_{\mathbf{n}}$ is maximally powerful in $\N^{\mathbf{n}}[S]$.
\end{theorem}

\noindent\emph{Proof}: This is provable as a meta-theorem of $EA_1$ with $\mathbf{n}$ as a meta-variable ranging over the classical natural numbers. We proceed by induction on $\mathbf{n}$.  The base case is trivial, and, to establish the inductive step, it suffices to show that, if an arithmetical global function, $\varphi(a)$, is maximally powerful in a natural number system, $\N$, then $S^{\varphi(\mathrm{log}_S(a))}$ is maximally powerful in $\N[S]$

First, note that, by a straightforward induction, if $\eta_{\varphi}(a) =_{\mathrm{df.}}  S^{\varphi(\mathrm{log}_S(a))}$,
\[
EA_1 \vdash \eta^{\mathbf{n}}_{\varphi}(a) \simeq_\mathrm{c} S^{\varphi^{\mathbf{n}}(\mathrm{log}_S(a))}
\]

\noindent Now, suppose $\varphi$ is maximally powerful in $\N$ and $\psi$ is an arithmetical global function such that $\N[S] \stackrel{\psi}{\rightarrow} \N[S]$.  Then define
\[
\psi'(a) =_{\mathrm{df.}}  \mathrm{max}_\mathrm{c}(\{\psi(x) : x \leq_\mathrm{c}  a\})
\]
Crucially,
\begin{enumerate}
\item[(i)] $\psi'$ is monotone increasing,
\item[(ii)] $EA_1 \vdash \psi(a) \leq_\mathrm{c} \psi'(a)$, and
\item[(iii)] A natural number system will be closed under $\psi'$ if, and only if, it is closed under $\psi$.
\end{enumerate}
By (ii), it will suffice to show that there is a number $\mathbf{n}$ such that
\[
\lambda x\, \psi'(x) \ll \lambda x\, S^{\varphi^{\mathbf{n}}(\mathrm{log}_S(x))}
\]

Now, by (iii), since $\N[S] \stackrel{\psi}{\rightarrow} \N[S]$, it follows that $\N[S] \stackrel{\psi'}{\rightarrow} \N[S]$, and thus, by Corollary \ref{cor2NNS}, $\N \stackrel{\lambda x\, \mathrm{log}_S(\psi'(S^x))}{\rightarrow} \N$.  Thus, since $\varphi$ is maximally powerful in $\N$, there is a number $\mathbf{n}$ such that $\lambda x\, \mathrm{log}_S(\psi'(S^x)) \ll \lambda x\, \varphi^{\mathbf{n}}(x)$.  Thus, there is a number $\mathbf{M}$ such that
\begin{enumerate}
\item[] $EA_1 \vdash\forall x[ \mathbf{M} \leq_\mathrm{c}  x \supset \mathrm{log}_S\left (\psi'\left(S^x\right )\right ) <_\mathrm{c}  \varphi^{\mathbf{n}}\left (x\right )]$, so
\item[] $EA_1\vdash \forall x[S^\mathbf{M} \leq_\mathrm{c}  x \supset \mathrm{log}_S\left (\psi'\left (S^{\mathrm{log}_S\left (x\right )}\right )\right ) <_\mathrm{c}  \varphi^{\mathbf{n}}\left (\mathrm{log}_S\left (x\right )\right )]$
\end{enumerate}
But, since $EA_1 \vdash \forall x[x \leq_\mathrm{c} S^{\mathrm{log}_S(x)}]$, by (i) we have
\begin{enumerate}
\item[] $EA_1 \vdash\forall x[ \mathrm{log}_S\left (\psi'\left (x\right )\right ) \leq_\mathrm{c}  \mathrm{log}_S\left (\psi'\left (S^{\mathrm{log}_S\left (x\right )}\right )\right )]$, so
\item[] $EA_1 \vdash \forall x[S^\mathbf{M} \leq_\mathrm{c}  x \supset \mathrm{log}_S\left (\psi'\left (x\right )\right ) <_\mathrm{c}  \varphi^{\mathbf{n}} \left (\mathrm{log}_S\left (x\right )\right )]$ . Thus
\item[] $EA_1 \vdash \forall x[S^\mathbf{M} \leq_\mathrm{c}  x \supset S^{\mathrm{log}_S\left (\psi'\left (x\right )\right )} <_\mathrm{c}  S^{\varphi^{\mathbf{n}}\left (\mathrm{log}_S\left (x\right )\right )}]$.  But
\item[] $EA_1 \vdash\forall x[ \psi'\left (x\right ) \leq_\mathrm{c}  S^{\mathrm{log}_S\left (\psi'\left (x\right )\right )}]$, so
\item[] $EA_1 \vdash\forall x[ S^\mathbf{M} \leq_\mathrm{c}  x \supset \psi'\left (x\right ) <_\mathrm{c}  S^{\varphi^{\mathbf{n}}\left (\mathrm{log}_S\left (x\right )\right )}]$, as required. \hfill $\Box$ \medskip
\end{enumerate}
Now, suppose we let $\mu_{\mathbf{n}}(a) =_{\mathrm{df.}} (\lambda x\, x +_\mathrm{c} \mathbf{1})^\mathbf{2}_{\mathbf{n}}(a)$.  Then, by the theorem just proved, if $\N$ is $\CH$, $\Z$, or $\VN$, then $\mu_{\mathbf{n}}$ is maximally powerful in $\N^{\mathbf{n}}[S]$, for $\lambda x\, x +_\mathrm{c} \mathbf{1}$ is maximally powerful in $\N$.  By an external induction along the classical natural numbers, we can show that, for any set $a$,
\[
\mu_{\mathbf{n}}(a) \simeq_\mathrm{c} \mathbf{2}^{\mathrm{log}^{\mathbf{n}}(a) + \mathbf{1}}_{\mathbf{n}}
\]

Using this representation of $\mu_\mathbf{n}$, it is straightforward to establish the following relationships between the first three functions in the $\mu_{\mathbf{n}}$ hierarchy and more familiar functions:
\begin{enumerate}
\item[(i)] $\lambda x\, \mathbf{2}x \ll \mu_\mathbf{1} \ll \lambda x\, \mathbf{4}x$
\item[(ii)] $\lambda x\, x^\mathbf{2} \ll \mu_\mathbf{2} \ll \lambda x\, x^\mathbf{4}$
\item[(iii)] $\lambda x\, \omega_1(x) \ll \mu_\mathbf{3} \ll \lambda x\, \omega^\mathbf{2}_1(x)$
\end{enumerate}
where, by definition, $\omega_1(a) =_\mathrm{df.} a^{\mathrm{log}_\mathrm{2}(a)}$.  This function was defined first by Wilkie and Paris \cite{awjp}.  The weak arithmetic, $I\Delta_0 + \Omega_1$, is obtained from $I \Delta_0$ by adding the formula that states that $\omega_1$ is total.  This theory has been studied extensively, in particular, with respect to its strong connections with the theory of complexity classes.  From the theorem just proved along with (i), (ii), and (iii), it follows that, if $\N$ is $\CH$, $\Z$, or $\VN$, then
\begin{enumerate}
\item[(i')] $\lambda x\, \mathbf{2}x$ is maximally powerful in $\N[S]$,
\item[(ii')] $\lambda x\, x^\mathbf{2}$ is maximally powerful in $\N^\mathbf{2}[S]$, and
\item[(iii')] $\lambda x\, \omega_1(x)$ is maximally powerful in $\N^\mathbf{3}[S]$.
\end{enumerate}
This concludes our description of the closure properties of natural number systems in a hierarchy of $S$-ary numerals.

The relations of length that hold between them are given by the following theorem, which was first proved in \cite{jm}, along with the observation that, for any natural number system $\N$ and any numeration base $S$, since $\N \stackrel{\lambda x\, S^x}{\rightarrow} \N[S]$, it is certainly the case that $\N \preceq \N[S]$.
\begin{theorem}
Suppose $\N$ is a natural number system and $S$ is a numeration base.  Then the following three propositions are equivalent:
\begin{enumerate}
\item[\emph{(1)}] $\N \stackrel{\lambda x\, S^x}{\rightarrow} \N$.
\item[\emph{(2)}] $\N \cong \N[S]$.
\item[\emph{(3)}] $\N[S] \stackrel{\lambda x\, S^x}{\rightarrow} \N[S]$.
\end{enumerate}
\end{theorem}
\emph{Proof}:  (2) follows from (1), (3) from (2), and (1) from (3) by Lemma \ref{NNS} and composing arrows. \hfill $\Box$\medskip

Thus, if $\N$ is not closed under $\lambda x\, S^x$, then, for all $\mathbf{n}$, $\N^{\mathbf{k}}[S]$ is not closed under $\lambda x\, S^x$ and
\[
\N \prec \N[S] \prec \N^{\mathbf{2}}[S] \prec \N^{\mathbf{3}}[S] \prec \cdots
\] 
On the other hand, if $\N$ is closed under $\lambda x\, S^x$,
\[
\N \cong \N[S] \cong \N^{\mathbf{2}}[S] \cong \N^{\mathbf{3}}[S] \cong \cdots
\]

\subsection{Fixed length numeral hierarchies}\label{nthrootext}

In this section, we describe a second hierarchy called the \emph{length $L$ numeral extension hierarchy}. Given a fixed linear ordering $L$ with two or more terms and a natural number system $\N$, we define its \emph{fixed length $L$ numeral extension}, $\mathcal{N}\langle L\rangle$, roughly as follows:  we fix the lengths of the numerals to be $L$ and allow the numeration base of the length $L$ numerals to range over $\{\mathrm{InSeg}(N) : N \mbox{ in } \N\mbox{ and }N\geq2\}$.\footnote{Recall that $\mathrm{InSeg}(N)$ is the set of all proper initial segments of $N$:  thus $\mathrm{InSeg}(L') \simeq_\mathrm{c} L'$.  As noted above, for any linear ordering with two or more terms, $L$, $\mathrm{InSeg}(L)$ is a numeration base.}

To define $\mathcal{N}\langle L\rangle$ we must specify both its initial term and its successor function. The initial term will be
\[
\left (\{[\, ], [0_\N]\}, \<[\, ], [\, ], ..., [\, ]\>_L \right)
\]
As is the case with all simply infinite systems the successor function for $\mathcal{N}\langle L\rangle$ operates on its terms and so can be regarded as enumerating those terms in their natural order. We shall arrange matters so that each term of $\mathcal{N}\langle L\rangle$ is a numeral of the form $(\mathrm{InSeg}(N),s)$, where $N\geq2$ lies in $\mathcal{N}$ and $s$ is of length $L$. Notice that our initial term for $\mathcal{N}\langle L\rangle$ is of this form.

The idea behind the construction of the sucessor function $\sigma_{\mathcal{N}\langle L\rangle}$ is to enumerate numerals to a given base of the form $S=\mathrm{InSeg}(N)$ in their natural order as long as we can until we reach the numeral
\[
(S,\<\max(S),\cdots,\max(S)\>_L)
\]
whose successor numeral to the base $S$ requires more than $L$ digits. We then continue the enumeration using numerals to the next largest base of the appropriate form, beginning the numeral to that base representing the next cardinality greater than the cardinality represented by the numeral $(S,\<\max(S),\cdots,\max(S)\>_L)$ (i.e, the numeral representing the cardinality of $S^L$ to the new base).

The definition of the natural number system $\mathcal{N}\langle L\rangle$ can now be given as follows:

\begin{definition}Let $\mathcal{N}$ be a natural number system and $L$ be a linear ordering of length at least two. Then $\mathcal{N}\langle L\rangle$ is generated from the initial element
\[
\left (\{[\, ], [0_\N]\}, \<[\, ], [\, ], ..., [\, ]\>_L \right)
\]
by the successor function $\sigma_{\mathcal{N}\langle L\rangle}$ defined by stipulating that for any numeral of the form $(S,s)$ where $S=\mathrm{InSeg}(N)$, for $N\in \mathcal{N}$:
\begin{enumerate}
\item[\emph{(1)}] If $s=\<d_0,\cdots,d_k\>_L$ where $d_k\neq\max(S)$, then
\[
\sigma_{\mathcal{N}\langle L\rangle}((S,s))=(\,S,\,\<d_0,\cdots,(d_k+1)\>_L\,)
\]
\item[\emph{(2)}] If $s=\<\max(S),\cdots,\max(S)\>_L$, so that $(S,s)$ is the largest base $S$ numeral of length $L$, then
\[
\sigma_{\mathcal{N}\langle L\rangle}((S,s))=(S',s')
\]
where
\begin{enumerate}
\item[\emph{(i)}]$S'=\mathrm{InSeg}(N')$, where $N'$ is the number of $\mathcal{N}$ next after $N$.
\item[\emph{(ii)}] $(S', s')$ is the unique length $L$ numeral to the base $S'$ that codes the cardinality of $S^L$.
\end{enumerate}
\end{enumerate}
\end{definition}
$\mathcal{N}\langle L\rangle$ so defined is clearly a natural number system. Moreover, for every $N$ lying in $\mathcal{N}$, the term $(S,\<\max(S),\cdots,\max(S)\>_L)$, where $S=\mathrm{InSeg}(N)$, lies in $\mathcal{N}\langle L\rangle$. Thus the cardinalities of the bases of the numerals that are the terms of $\mathcal{N}\langle L\rangle$ range over all cardinalities of numbers $N\geq2$ in $\mathcal{N}$.

As before, we define, by meta-recursion, the following meta-sequence of natural number systems:
\[
\N^{\mathbf{0}}\langle L\rangle= _\mathrm{df.}  \N\,\,\,\,\,\,\,\,\,\,\,\,\,\,\,   \N^{\mathbf{k+1}}\langle L\rangle =_\mathrm{df.}  (\N^{\mathbf{k}}\langle L\rangle)\langle L\rangle
\]
and we call the resulting hierarchy the \emph{hierarchy of length $L$ numeral extensions of $\N$}.

For each result concerning $S$-ary numeral extensions that we proved above, there is an analogous fact concerning length $L$ numeral extensions.  We list these facts below, but omit proofs except in a few special cases.

The next Lemma concerning length $L$ numeral extensions is analogous to Lemma \ref{NNS} concerning base $S$ numeral extensions:
\begin{lemma}\label{multext} $\N \stackrel{\lambda x\, x^L}{\rightarrow} \mathcal{N}\langle L\rangle$
and $\mathcal{N}\langle L\rangle \stackrel{\lambda x\, \sqrt[L]{x}}{\rightarrow} \N$
\end{lemma}
\emph{Proof}:  To show that $\N \stackrel{\lambda x\, x^L}{\rightarrow} \mathcal{N}\langle L\rangle$, let $N$ lie in $\mathcal{N}$ and let $S=\mathrm{InSeg}(N)$. Then $[\cdots,(S,\<\max(S),\cdots\max(S)\>_L)]_{\mathcal{N}\langle L\rangle}$ lies in $\mathcal{N}\langle L\rangle$ and has $S^L$ members. So
\[
N^L\simeq_\mathrm{c}S^L\simeq_\mathrm{c}[\cdots,(S,\<\max(S),\cdots\max(S)\>_L)]_{\mathcal{N}\langle L\rangle}
\]
\medskip

\noindent Next, we show that $\mathcal{N}\langle L\rangle \stackrel{\lambda x\, \sqrt[L]{x}}{\rightarrow} \N$.  Suppose $N = [0, ..., k]$ lies in $\N$ and $[\cdots, \left(\mathrm{InSeg}([0, ..., k]), s \right )]_{\N \langle L \rangle}$ is in $\N \langle L \rangle$.  Then
\begin{eqnarray*}
(N-1)^L & \simeq_\mathrm{c} & \left [\cdots, \left(\mathrm{InSeg}([0, ..., k-1]), \<\max(S), ..., \max(S)\>_L \right ) \right]_{\N[L]} \\
& \subsetneqq_* & \left [\cdots, \left(\mathrm{InSeg}([0, ..., k]), s \right )\right ]_{\N[L]} \\
& \subseteq_* & \left [\cdots, \left(\mathrm{InSeg}([0, ..., k]), \<\max(S), ..., \max(S)\>_L \right ) \right]_{\N[L]} \\
& \simeq_\mathrm{c} & N^L
\end{eqnarray*}
Thus,
\[
\sqrt[L]{ \left [\cdots, \left(\mathrm{InSeg}(N), s \right )\right ]} \simeq_\mathrm{c} N
\]
So let $\eta_2\left ( \left [\cdots, \left(\mathrm{InSeg}(N), s \right )\right ] \right ) = N$.  Then $\eta_2 : \mathcal{N}\langle L\rangle \stackrel{\lambda x\, \sqrt[L]{x}}{\rightarrow} \N$.  \hfill $\Box$

\begin{lemma}\label{lemNtimesmaxpow}
Suppose $\N$ is a natural number system and $\varphi$ an arithmetical global function.  Then $\varphi$ is maximally powerful in $\N$  if, and only if, $\lambda x\, (\varphi(\sqrt[L]{x}))^{L}$ is maximally powerful in $\N\langle L \rangle$.
\end{lemma}
Thus, we define the following hierarchy of arithmetical global functions.
\begin{definition}
Suppose $\varphi$ is an arithmetical global function; then
\[
\begin{array}{lll}
\varphi^L_\mathbf{0}(a) & = & \varphi(a)\\
\varphi^L_\mathbf{k+1}(a) & = & (\varphi_\mathbf{k}(\sqrt[L]{a}))^{L}
\end{array}
\]
\end{definition}
\begin{corollary}\label{corvarnNtimes}
Suppose $\N$ is a natural number system in which the arithmetical global function $\varphi$ is maximally powerful.  Then, for all $\mathbf{k}$ in $\mathbb{N}$, $\varphi^L_\mathbf{k}(x)$ is maximally powerful in $\N^{\mathbf{k}}\langle L \rangle$.
\end{corollary}
Now, suppose we define a hierarchy of functions $\xi_{\mathbf{n}, \mathbf{k}}$ as follows:
\begin{eqnarray*}
\xi^L_\mathbf{0}(a) & = & a + 1 \\
\xi^L_\mathbf{k+1}(a) & = & \left (\xi^L_\mathbf{k}\left (\sqrt[L]{a}\right)\right )^{L}
\end{eqnarray*}
Then, if $\lambda x\, x + \mathbf{1}$ is maximally powerful in $\N$, then $\xi^L_\mathbf{k}$ is maximally powerful in $\N^\mathbf{k}\langle L \rangle$.  By an (external) induction along the classical natural numbers, for any set $a$,
\[
\xi^L_\mathbf{k}(a) = \left (\sqrt[L^\mathbf{k}]{a} + \mathbf{1}\right )^{L^\mathbf{k}}
\]

This completes our description of the closure properties of the systems in the $L$-long numeral extensions hierarchy.  Thus, we turn to the relations of length between those systems.  In fact, our investigation will bring us again to consider the closure properties of those systems, as we shall see.

Again, in analogy with the case of $\N[S]$, we have $\N \preceq \mathcal{N}\langle L\rangle$, as well as the following result:

\begin{theorem}\label{mayberry3mult}  Suppose $\N$ is a natural number system.  Then the following three statements are equivalent:
\begin{enumerate}
\item[\emph{(1)}] $\mathcal{N}\langle L\rangle$ is closed under $\lambda x\, x^{\mathbf{2}}$.
\item[\emph{(2)}] $\mathcal{N}\langle L\rangle \preceq \N$ \emph{(}i.e., $\N$ measures $\mathcal{N}\langle L\rangle$\emph{)}.
\item[\emph{(3)}] $\N$ is closed under $\lambda x\, x^{\mathbf{2}}$.
\end{enumerate}
\end{theorem}
Thus, if $\N$ is not closed under $\lambda x\, x^{\mathbf{2}}$, then, for all natural numbers $\mathbf{k}$, $\N^{\mathbf{k}}\langle L \rangle$ is not closed under $\lambda x\, x^{\mathbf{2}}$ and
\[
\N \prec \mathcal{N}\langle L\rangle \prec \N^{\mathbf{2}}\langle L \rangle \prec \N^{\mathbf{3}}\langle L \rangle \prec \cdots
\] 
On the other hand, if $\N$ is closed under $\lambda x\, x^{\mathbf{2}}$, we have
\[
\N \cong \mathcal{N}\langle L\rangle \cong \N^{\mathbf{2}}\langle L \rangle \cong \N^{\mathbf{3}}\langle L \rangle \cong \cdots
\]
Thus, just as we cannot obtain a system closed under $\lambda x\, S^x$ by the technique of $S$-ary numeral extension unless we have one already, so we cannot obtain a natural number system closed under $\lambda x\, x^{\mathbf{2}}$ by the technique of $L$-long numeral extensions unless we have one already.

However, there is a curious additional fact about $L$-long numerals that is not analogous to any known fact about $S$-ary numerals.  This theorem and Theorem \ref{boundedbyaddition} concern $L$-long numeral extensions where $L \simeq_\mathrm{c} \mathbf{n}$ for some metatheoretical $\mathbf{n}$.
\begin{theorem}
Suppose $\N$ is a natural number system, and suppose $L \simeq_\mathrm{c} \mathbf{n}$.  Then the following two statements are equivalent:
\begin{enumerate}\label{mayberry4mult}
\item [\emph{(1)}] $\mathcal{N}\langle L\rangle$ is closed under $\lambda x\, \mathbf{2}x$.
\item [\emph{(2)}] $\N$ is closed under $\lambda x\, \mathbf{2}x$.
\end{enumerate}
\end{theorem}
\emph{Proof. }  Suppose (1).  Then, since $L \simeq_\mathrm{c} \mathbf{n}$, we have $\mathcal{N}\langle L\rangle \stackrel{\lambda x\, \mathbf{2}^{L}}{\rightarrow} \mathcal{N}\langle L\rangle$ and
\[
\N \stackrel{\lambda x\, x^{L}}{\rightarrow} \mathcal{N}\langle L\rangle \stackrel{\lambda x\, \mathbf{2}^{L}x}{\rightarrow} \mathcal{N}\langle L\rangle \stackrel{\lambda x\, \sqrt[L]{x}}{\rightarrow} \N
\]
Thus, $\N \stackrel{\lambda x\, \sqrt[L]{\mathbf{2}^{L}x^{L}}}{\rightarrow} \N$, which entails (2) since $\sqrt[L]{\mathbf{2}^{L}x^{L}} = \mathbf{2}x$.

\medskip

\noindent Suppose (2).  Then
\[
\mathcal{N}\langle L\rangle \stackrel{\lambda x\, \sqrt[L]{x}}{\rightarrow} \N \stackrel{\lambda x\, \mathbf{2}x}{\rightarrow} \N \stackrel{\lambda x\, x^{L}}{\rightarrow} \mathcal{N}\langle L\rangle
\]
Thus, $\mathcal{N}\langle L\rangle \stackrel{\lambda x\, \left (\mathbf{2}\sqrt[L]{x}\right )^{L}}{\rightarrow} \mathcal{N}\langle L\rangle$, so $\mathcal{N}\langle L\rangle$ is closed under $\lambda x\, \mathbf{2}^{L}x$ and thus under $\lambda x\, \mathbf{2}x$. \hfill $\Box$

\medskip

Thus, if $L \simeq_\mathrm{c} \mathbf{n}$ and if $\N$ is not closed under $\lambda x\, \mathbf{2}x$, then for all natural numbers $\mathbf{k}$, $\N^{\mathbf{k}}\langle L \rangle$ is not closed under $\lambda x\, \mathbf{2}x$.  Thus, not only can we not obtain a natural number system closed under $\lambda x\, x^{\mathbf{2}}$ by the technique of $L$-long numeral extensions unless we have one already; if we are allowed only numerals of classically finite length, neither can we thus obtain a natural number system closed under $\lambda x\, \mathbf{2}x$ unless we have one already.

This result reveals the possibility of taking a natural number system, $\N$, in which $\lambda x\, \mathbf{2}x$ is maximally powerful (e.g. $\Z[\mathbf{2}]$) and which is thus not closed under $\lambda x\, x^{\mathbf{2}}$ and then taking successive length $L$ numeral extensions to give the infinite hierarchy,
\[
\N \prec \mathcal{N}\langle L\rangle \prec \N^{\mathbf{2}} \langle L \rangle \prec \N^{\mathbf{3}} \langle L \rangle \prec \cdots
\]
of natural number systems that is guaranteed by Theorem \ref{mayberry3mult}.

Thus, the commensurability relations in this hierarchy are exactly the same as when, above, we built an infinite hierarchy on a natural number system $\N$ such that $\lambda x\, x + \mathbf{1}$ is maximally powerful in $\N$:  that is, each of the natural number systems in this new hierarchy is strictly longer than its predecessors. However, in the case where $\lambda x\, x + \mathbf{1}$ was maximally powerful in $\N$, each $\N^{\mathbf{n}}\langle L \rangle$ was closed under $\xi_{\mathbf{n}}$ and each $\xi_{\mathbf{n}}$ was genuinely more powerful than $\xi_{\mathbf{n-1}}$.  Thus, if $\lambda x\, x + \mathbf{1}$ is maximally powerful in $\N$, each $\N^{\mathbf{n}}\langle L \rangle$ is closed under new, stronger functions than $\N^{\mathbf{n-1}}\langle L \rangle$.  On the other hand, when we take $\N$ such that $\lambda x\, \mathbf{2}x$ is maximally powerful in $\N$, and $L$ is classically finite, the $L$-long numeral extensions, $\mathcal{N}\langle L\rangle$, $\N^{\mathbf{2}}\langle L \rangle$, etc., are not closed under any new arithmetical global functions.  This is the content of the following result.
\begin{theorem}\label{boundedbyaddition}
Suppose $\lambda x\, \mathbf{2}x$ is maximally powerful in the natural number system $\N$ and suppose $L \simeq_\mathrm{c} \mathbf{n}$ for metatheoretical $\mathbf{n}$.  Then, for all natural numbers $\mathbf{k}$, if $\N^{\mathbf{k}}\langle L \rangle$ is closed under the arithmetical global function $\varphi$, then there is $\mathbf{K}$ in $\mathbb{N}$ such that
\[
EA_1 \vdash \forall x[\varphi(x) < \mathbf{K}x]
\]
\end{theorem}
\emph{Proof. } We proceed by (meta)-induction on $\mathbf{k}$:

\medskip

\noindent \textbf{Basis Case ($\mathbf{k}=\mathbf{0}$) } By definition, since $\lambda x\, \mathbf{2}x$ is maximally powerful in $\N$, the hypothesis is true for $\mathbf{k} = \mathbf{0}$.

\medskip

\noindent \textbf{Inductive Step ($\mathbf{k} > \mathbf{0}$)} Suppose that, whenever $\N^{\mathbf{k}}\langle L \rangle$ is closed under $\varphi$, there is a natural number $\mathbf{K}$ such that
\[
EA_1 \vdash \forall x[\varphi(x) <_\mathrm{c}  \mathbf{K}x]
\]
And suppose $\N^{\mathbf{k+1}}\langle L \rangle$ is closed under $\psi$.  First, define the arithmetical global function $\psi'$ as follows:
\[
\psi'(a) = \mathrm{max}_\mathrm{c}(\{\psi(x) : x \leq_\mathrm{c} a\})
\]
and then define the arithmetical global function $\psi''$ as follows:
\[
\psi''(a) = \psi'\left ( \left (\sqrt[L]{a} + \mathbf{1}\right )^{L} \right )
\]
Then the important facts about $\psi''$ are these:
\begin{enumerate}
\item[(1)] $EA_1 \vdash\forall x[ \psi(x) \leq_\mathrm{c} \psi''(x)]$;
\item[(2)] $EA_1 \vdash\forall x[ \psi''(x) \leq_\mathrm{c}  \psi''((\sqrt[L]{a} + \mathbf{1})^{L})]$
\item[(3)] If $\M$ is a natural number system closed under $\lambda x\, \mathbf{2}x$, then $\M$ is closed under $\psi''$ if, and only if, it is closed under $\psi$.
\end{enumerate} 
(1) and (2) are obvious.  It is also obvious that, for any natural number system $\M$, $\M$ is closed under $\psi'$ if, and only if, $\M$ is closed under $\psi$.  By (1), we have that, if $\M$ is closed under $\psi''$ it is closed under $\psi$.  This gives one half of (3).  For the converse, suppose $\M$ is closed under $\lambda x\, \mathbf{2}x$ and under $\psi$ and suppose $a$ is a segment of $\M$.  Then, since 
\[
a \leq_\mathrm{c} (\sqrt[L]{a} + \mathbf{1})^{L} \leq_\mathrm{c} L.L!a \simeq_\mathrm{c} (\mathbf{n +1})(\mathbf{n +1})!a,
\]
$\M$ is closed under $\lambda x\, (\sqrt[L]{a} + \mathbf{1})^{L}$ and so, by composition, $\M$ is closed under $\psi''$.  Thus, (3).

Then
\[
\N^{\mathbf{k}}\langle L \rangle \stackrel{\lambda x\, x^{L}}{\rightarrow} \N^{\mathbf{k+1}}\langle L \rangle \stackrel{\psi''}{\rightarrow} \N^{\mathbf{k+1}}\langle L \rangle \stackrel{\lambda x\, \sqrt[L]{x}}{\rightarrow} \N^{\mathbf{k}}\langle L \rangle
\]
Thus, by composition, $\N^{\mathbf{k}}\langle L \rangle$ is closed under $\lambda x\, \sqrt[L]{\psi''(x^{L})}$.  Thus, by inductive hypothesis, we have that there is natural number $\mathbf{K}$ such that
\begin{enumerate}
\item[] $EA_1 \vdash\forall x[ \sqrt[L]{\psi''(x^{L})} <_\mathrm{c}  \mathbf{K}x]$, so
\item[] $EA_1 \vdash\forall x[ \sqrt[L]{\psi''(x^{L})} + \mathbf{1} <_\mathrm{c}  (\mathbf{K+1})x]$.  Thus,
\item[] $EA_1 \vdash\forall x[ \psi''(x^{\mathbf{n}}) <_\mathrm{c} (\sqrt[L]{\psi''(x^{L})} + \mathbf{1})^{L} <_\mathrm{c}  (\mathbf{K+1})^{L}x^{L}]$.  But then
\item[] $EA_1 \vdash\forall x[ \psi(x) \leq_\mathrm{c} \psi''(x) \leq_\mathrm{c}  \psi''(( \sqrt[L]{x} + \mathbf{1})^{L})$

\hspace{9mm} $ <_\mathrm{c} (\mathbf{K+1})^{L}( \sqrt[L]{x} + \mathbf{1})^{L} \leq_\mathrm{c} (\mathbf{K+1})^{L}L.L!x$

\hspace{9mm} $\simeq_\mathrm{c} (\mathbf{K+1})^\mathbf{n+1}(\mathbf{n+1}).(\mathbf{n+1})!x]$
\end{enumerate}
and our induction is complete. \hfill $\Box$

\begin{theorem}
If $\N$ is closed under $\lambda x\, \mathbf{2}x$ and $L \simeq_\mathrm{c} \mathbf{n}$, then, for all natural numbers $\mathbf{k}$, $\N^{\mathbf{k}}\langle L \rangle$ and $\N$ are closed under the same arithmetical global functions.
\end{theorem}
\emph{Proof. } If $\N^{\mathbf{k}}\langle L \rangle$ is closed under $\varphi$, then there is natural number $\mathbf{K}$ such that
\[
EA_1 \vdash \forall x[\varphi(x) \leq_\mathrm{c} \mathbf{K}x]
\]
and thus, since $\N$ is closed under $\lambda x\, \mathbf{2}x$ and thus under $\lambda x\, \mathbf{K}x$, $\N$ is closed under $\varphi$ as required. \hfill $\Box$\medskip

\section{Lexicographic ordering}\label{lexordsect}

As we noted in Section \ref{EAweakarith}, in \cite{wa}, Ackermann gave a coding of the hereditarily finite pure sets in the natural numbers.  This coding induces an ordering on the hereditarily finite pure sets via the standard ordering of their Ackermann codes as natural numbers.

We intend to study this ordering, but since we don't have natural numbers to function as codes, we must define the ordering on sets determined by their Ackermann codes directly in set theory using the concept of \emph{lexicographic orderings}.
\begin{definition}\label{lexorddef} Let $L$ be a linear ordering. Then $\Lex(L)$ is the linear ordering of $\P(\mathrm{Field}(L))$ given by stipulating that for $X,Y \in \mathcal{P}(\mathrm{Field}(L))$
\[
X\leq_{\Lex(L)}Y\Leftrightarrow_{\mathrm{df.}} \mbox{the $L$-greatest element of $X \bigtriangleup Y$ is in $Y$}
\]
where $X \bigtriangleup Y =_{\mathrm{df.}} (X - Y) \cup (Y-X)$ is the symmetric difference function.
\end{definition} 
Thus, if $L = [\varnothing, \{\varnothing\}, \{\{\varnothing\}\}]$, then $X = \{\varnothing, \{\varnothing\}\} \leq_{\Lex(L)} \{\{\varnothing\}, \{\{\varnothing\}\}\} = Y$ since $X \bigtriangleup Y = \{\varnothing, \{\{\varnothing\}\}\}$, $\varnothing <_L \{\{\varnothing\}\}$, and $\{\{\varnothing\}\} \in Y$.\medskip

Lexicographic orderings have the following important property:
\begin{proposition}\label{ACKfund}Let $L$ and $L'$ be linear orderings with $ L \subsetneqq_* L'$. Then
\[
\Lex(L) \subsetneqq_* \Lex(L')
\]  
\end{proposition}

\begin{definition} The \emph{lexicographic hierarchy system}, $\LEX$,  is the iteration system whose initial term is the empty ordering $[\,]$ and whose iterating function $\sigma_\LEX$ is $\Lex$.
\end{definition}
By Proposition \ref{ACKfund}, $\LEX$ is a natural number system.  Moreover, every term of $\LEX$ is a linear ordering.

\begin{proposition} $\CH \cong \LEX$
\end{proposition}
\emph{Proof}:  First, we show that $\LEX \preceq \CH$.  For any linear ordering $L$,
\[
\Field(\Lex(L)) = \P(\Field(L))
\] 
Thus, if $[L_0,\,L_1,\,\cdots,\,L_k]$ is a number in $\LEX$, then each $\Field(L_i)$ is a term of $\CH$ and $[\mathrm{Field}(L_0),\,\mathrm{Field}(L_1),\,\cdots,\,\mathrm{Field}(L_k)]$ is a number of  $\CH$.

Next, we show that $\CH \preceq \LEX$. Given a $\CH$-number $[\varnothing = V_0,\,\cdots,\,V_n]$, we can use recursion along $[V_0,\,\cdots,\,V_n]$ to define $[[\, ] = L_0, \cdots, L_n]$, a $\LEX$-number of the same length as $[V_0,\,\cdots,\,V_n]$.   Let $S$ be the set of linear orderings of subsets of $V_n$ and define $g : \{V_0, \cdots, V_n\} \rightarrow S$ as follows:
\begin{eqnarray*}
g ` V_0 & = & [\, ] \\
g `V_{i+1} & = & \Lex(g `V_i)
\end{eqnarray*}
Then let $L_i = g ` V_i$.  Then $[[\, ] = L_0, \cdots, L_n]$ is a number in $\LEX$ and $[[\, ] = L_0, \cdots, L_n] \simeq_\mathrm{c} [\varnothing = V_0, \cdots, V_n]$ as required. \hfill $\Box$\medskip

Thus the natural number systems $\LEX$ and $\CH$ are of the same length.  What's more, given a $\CH$-number $[V_0,\cdots,V_n]_\CH$ and its corresponding $\LEX$-number, $[L_0,\cdots,L_n]_\LEX$, for $i = 0, \cdots, n$, $\Field(L_i) = V_i$.\footnote{Although we have used the `indices' $i = 0, \cdots, n$ to describe the $\CH$-number $[V_0,\cdots,V_n]_\CH$ and its corresponding $\LEX$-number, $[L_0,\cdots,L_n]_\LEX$, this is just a device to simplify this remark. In fact the remark could be made `inside' $EA_1$, so to speak.}

Now we will use $\LEX$ to define a further natural number system, $\ACK$, whose terms are the initial segments of the terms of $\LEX$.

\begin{definition} The \emph{Ackermann  system} $\ACK$ is the iteration system whose initial term is $[\, ]\,(=\varnothing)$, and whose iterating function is defined as follows:  If $l$ is a linear ordering, and if there is a term $L$ of $\LEX$ such that $L \subseteq_* l \subsetneqq_* \Lex(L)$, then let 
\[
\sigma_{\ACK}(l) =_{\mathrm{df.}}  \mbox{the unique } l' \subseteq_* \Lex(L) \mbox{ such that } \mathrm{Field}(l') \simeq_\mathrm{c} \mathrm{Field}(l) +_\mathrm{c} \mathbf{1}
\]
Otherwise, let $\sigma_\ACK(l) = \varnothing$; i.e. a `don't care' value.
\end{definition}
It is a straightforward consequence of Proposition \ref{ACKfund} that $\ACK$ is a natural number system.  Moreover, it is clear that the numbers of $\ACK$ are recoverable from its terms: given a term $l_{k+1}=[s_0,\,\cdots, s_k]$ of $\ACK$, then the linear ordering
\[
[\,\,\,[\, ]\, (\,= l_0),\,[s_0]\, (\,= l_1),\,[s_0,s_1]\,(\,= l_2),\,\cdots,\,[s_0,s_1,s_2,\cdots,s_k]\,(\, = l_{k+1})\,\,\,]
\]
of its initial segments is the $\ACK$-number whose last term is $l_{k+1}$.  Thus, $[l_0, \cdots, l_k] \simeq_\mathrm{c} l_k + \mathbf{1}$.  Furthermore, a straightforward induction establishes that the field of every $\ACK$-term is transitive.

There is a close connection between $\ACK$ and $\LEX$:  by a straightforward induction, we can show that every term of $\LEX$ is also a term of $\ACK$; of course, between each $\LEX$-term regarded as an $\ACK$-term, there lie many further $\ACK$-terms; for this reason, $\LEX \preceq \ACK$.

$\ACK$ seems somewhat complicated. The underlying `Dedekind sequence'
\[
l_0,\, l_1,\,l_2,\,\cdots,\,l_k,\,\cdots
\]
consisting of the \emph{terms} of $\ACK$ ordered globally in the manner determined by the \emph{numbers} of $\ACK$ is itself the class of initial segments of the underlying `Dedekind sequence'
\[
s_0\,(=\varnothing),\,s_1\,(=\{\varnothing\}),s_2\,(=\{\{\varnothing\}\}),\,\cdots
\]
which lists finite pure sets in the order determined by their Ackermann codes.

It is thus tempting to regard the species of $\ACK$-\emph{terms}, arranged in the obvious order determined by $\ACK$, as the numbers of a natural number system $\ACK_0$, more basic than the system $\ACK$ by virtue of the fact that its underlying `Dedekind sequence' is the ordering of the hereditarily finite pure sets determined by their Ackermann codes. Under that description, however, $\ACK_0$ is not a proper natural number system: we have not specified it in terms of an initial term together with a successor function $\sigma_{\ACK_0}$.  However, we can specify it in this way by means of the following definition.

\begin{definition}
Let $\ACK_0$ be the iteration system generated from $\varnothing$ by the following iterating function:  Given a set $S$, let $[s_0, \cdots, s_k]$ be the longest term of $\ACK$ such that $\{s_0, \cdots, s_k\} \subseteq S$; then, by definition
\[
\sigma_{\ACK_0}(S) = (S - \{s_0, \cdots, s_k\}) \cup \{s_{k+1}\}
\]
where $s_{k+1} = \Last(\sigma_{\ACK}([s_0, \cdots, s_k]))$.
\end{definition}

To show that the numbers of $\ACK_0$ are the terms of $\ACK$, we must show that, if $[s_0, \cdots, s_n]$ is a term of $\ACK$, then 
\[
\sigma_{\ACK}([s_0, \cdots, s_n]) = [s_0, \cdots, s_n, \sigma_{\ACK_0}(s_n)]
\]
To do this, we must prove the following:

\begin{theorem}
Suppose $[s_0, \cdots, s_n]$ is a term of $\ACK$ and $L$ is the longest term of $\LEX$ such that $L \subseteq_* [s_0, \cdots, s_n] \subsetneqq_* \Lex(L)$. Then
\begin{enumerate}
\item[\emph{(i)}] $s_n <_{\Lex(L)} \sigma_{\ACK_0}(s_n)$; and
\item[\emph{(ii)}] If $a \subseteq \Field(L)$ and $s_n <_{\Lex(L)} a$, then $\sigma_{\ACK_0}(s_n) \leq_{\Lex(L)} a$.
\end{enumerate}
\end{theorem}
\emph{Proof}: (i) We must show that the $L$-greatest element of $s_n \bigtriangleup \sigma_{\ACK_0}(s_n)$ lies in $\sigma_{\ACK_0}(s_n)$.  Now, if $[s_0, \cdots, s_k]$ is the longest Ackermann ordering such that $\{s_0, \cdots, s_k\} \subseteq s_n$, then $s_n \bigtriangleup \sigma_{\ACK_0}(s_n) = \{s_0, \cdots, s_k, s_{k+1}\}$, and the $L$-greatest element in this set is $s_{k+1}$, which is in $\sigma_{\ACK_0}(s_n)$.\smallskip

\noindent(ii) Suppose $a \subseteq \Field(L)$ and $s_n <_{\Lex(L)} a$ and suppose again that $[s_0, \cdots, s_k]$ is the longest Ackermann ordering such that $\{s_0, \cdots, s_k\} \subseteq s_n$.  Then let $s_{m_0}$ be $L$-greatest element of $a \bigtriangleup s_n$.  Then, since $s_n <_{\Lex(L)} a$, $s_{m_0} \in a$.  Also, $s_{m_0}$ is not in $\{s_0, \cdots, s_k\} \subseteq s_n$:  if it were, then it would be in $a$ and in $s_n$ and thus not in their symmetric difference.  Thus, $s_0, \cdots, s_k <_L s_{m_0}$.  Thus, there are two cases:
\begin{enumerate}
\item[(a)] Suppose $s_{m_0} = s_{k+1}$.  Since $s_{k+1}$ is the $L$-greatest element of $a \bigtriangleup s_n$ and by the definition of $\sigma_{\ACK_0}$, if $s_{k+1} < s_p$, then
\[
s_p \in a \equiv s_p \in s_n
\]
If not, $s_p$ would be in their symmetric difference, and it is $L$-greater than $s_{m_0}$.
Also, if $s_{k+1} < s_p$, then, since $s_n$ and $\sigma_{\ACK_0}(s_n)$ agree on elements above $s_{k+1}$,
\[
s_p \in s_n \equiv s_p \in \sigma_{\ACK_0}(s_n).
\]
Thus, if $s_{k+1} < s_p$, then $s_p \in a$ if, and only if, $s_p \in \sigma_{\ACK_0}(s_n)$.
It follows that $a \bigtriangleup \sigma_{\ACK_0}(s_n)$ contains no elements that are $L$-greater than $s_k$.  But $\sigma_{\ACK_0}(s_n)$ contains no elements that are \emph{not} $L$-greater than $s_k$.  Therefore, $a \bigtriangleup \sigma_{\ACK_0}(s_n)\subseteq a$ and thus, the $L$-greatest element of $a \bigtriangleup \sigma_{\ACK_0}(s_n)$ lies in $a$, as required.
\item[(b)] Suppose $s_{k+1} <_L s_{m_0}$.  Then $s_{m_0}$ is the $L$-greatest element of $a \bigtriangleup \sigma_{\ACK_0}(s_n)$ since, if $s_{k+1} < s_p$, then
\[
s_p \in \sigma_{\ACK_0}(s_n) \equiv s_p \in s_n
\] 
\end{enumerate}
This completes our proof of (ii) and the theorem follows. \hfill $\Box$ \medskip

Thus, the \emph{numbers} of $\ACK_0$ are the terms of $\ACK$.  However, the resulting natural number system $\ACK_0$ suffers from the disadvantage that we cannot recover its \emph{numbers} from its \emph{terms}: that is, there is no global function $\psi$ such that, if $[s_0, \cdots, s_n]$ is a number of $\ACK_0$, then $\psi(s_n) = [s_0, \cdots, s_n]$.  In fact, it is not even possible to define an inverse to the successor function for $\ACK_0$ on $\ACK_0$:  that is we can establish the following theorem

\begin{theorem}\label{noACK0pred}
There is no global function $\psi$ such that for all $\ACK_0$-numbers $[s_0, \cdots, s_k]$, $\psi(s_k) = s_{k-1}$.
\end{theorem}
\emph{Proof}: We show that, if there were such a global function, $\psi$, there would be a global function, $\psi'$, such that for all numbers $\mathbf{m}$, $\psi'(\mathbf{m}_{\Z}) = \mathbf{m}_{\VN}$.  This contradicts Lemma \ref{VNZconstant}.

For any number $\mathbf{m}$, $\mathbf{m}_{\Z}$ is always a singleton.  So, by the definition of $\sigma_{\ACK_0}$, if $\psi$ is an inverse for $\sigma_{\ACK_0}$ on the terms of $\ACK_0$, then $\psi(\mathbf{m}_{\Z})$ is the field of a term $[s_0, \cdots, s_n]$ of $\ACK$, and $\mathbf{m}_{\Z} = \{s_{n+1}\}$.  Thus, if $L$ is the longest term of $\Lex$ such that $L \subseteq_* [s_0, \cdots, s_n] \subsetneqq_* \Lex(L)$, then, since $\Field(\Lex(L))$ is a term of the cumulative hierarchy system and $\mathbf{m}_{\Z} \in \Field(\Lex(L))$, it follows that $\mathbf{m}_{\VN} \in \Field(\Lex(L))$, as required. \hfill $\Box$\medskip

Thus, $\ACK_0$ is a natural number system for which the $\Sigma_1$ definition in $EA_1$ for being a term of that system cannot be replaced by a quantifier-free definition in $EA_0$ that exploits the existence of a global function that recovers numbers from terms.


Let us now turn to the questions of closure and commensurability that concern $\ACK$.

\begin{theorem}\label{ACKclosure} $\lambda x\, \mathbf{2}^x$ is maximally powerful in $\ACK$.
\end{theorem}
\noindent\emph{Proof}:  It will suffice to show that $\ACK$ is closed under $\lambda x\, \mathbf{2}^x$, for it is a consequence of the Global Function Bounding Lemma (in conjunction with the existence of a global function that takes a set to a transitive set of the same size) that, for every global function $\varphi$ of $EA_0$, there is $\mathbf{n}$ such that $\varphi \ll \lambda x\, \mathbf{2}^x_{\mathbf{n}}$.

Suppose $[l_0, \cdots, l_n]$ is a number of $\ACK$ and $L_k$ and $L_{k+1}$ are successive terms of $\LEX$ such that $L_k \subseteq_* l_n \subsetneqq_* L_{k+1}$.  Then $[\cdots, L_{k+\mathbf{3}}]_\ACK$ is a number of $\ACK$ and
\begin{eqnarray*}
[\cdots, L_{k+\mathbf{3}}]_\ACK & \simeq_\mathrm{c} & \Lex(\Lex(L_{k+\mathbf{1}})) + \mathbf{1} \simeq_\mathrm{c} \mathbf{2}^{\mathbf{2}^{L_{k+1}}} \\
& > & \mathbf{2}^{L_{k+1} + \mathbf{1}} \simeq_\mathrm{c} \mathbf{2}^{[\cdots, L_{k+1}]} >\mathbf{2}^{[l_0, \cdots, l_n]}
\end{eqnarray*}
as required.   \hfill $\Box$\medskip

Thus $\ACK$ is  closed under all arithmetical functions, and this provides a complete answer to questions of closure in $\ACK$. Now let us now consider the relationships in respect of length between $\ACK$ and other natural number systems that we have defined. First we need to establish the following result:

\begin{proposition}\label{CHACK} $\ACK\stackrel{\lambda x\, \mathrm{suplog}_\mathbf{2}(x)}{\rightarrow}\CH, \VN, \Z$.
\end{proposition}
\noindent\emph{Proof}:  Since $\CH \preceq \VN, \Z$, it suffices to show that $\ACK \stackrel{\lambda x\, \mathrm{suplog}_\mathbf{2}(x)}{\rightarrow}\CH$.  If $[l_0, \cdots, l_n]$ is a number of $\ACK$ and $L_k$ and $L_{k+1}$ are successive terms of $\LEX$ such that $L_k \subseteq_* l_n \subsetneqq_* L_{k+1}$, then $[V_0, \cdots, V_k]$ is a number of $\CH$ and
\[
\mathbf{2}_{[V_0, \cdots, V_k]} \simeq_\mathrm{c} V_{k+\mathbf{2}} \simeq_\mathrm{c} L_{k+ \mathbf{2}}> L_{k+\mathbf{1}} \geq l_n + \mathbf{1} \simeq_\mathrm{c} [l_0, \cdots, l_n] 
 \]
 Thus, $[V_0, \cdots, V_k] \geq \mathrm{suplog}_{\mathbf{2}}([l_0, \cdots, l_n])$, which entails (ii).
 \hfill $\Box$\medskip

Now we can establish the length relations among $\CH$,$\VN$, $\Z$, and $\ACK$. 
\begin{theorem}\ 
\begin{enumerate}
\item[\emph{(i)}] $\CH \preceq \ACK$
\item[\emph{(ii)}] $\ACK \not \preceq \CH, \VN,\Z$
\item[\emph{(iii)}] $\VN, \Z \not \preceq \ACK$
\end{enumerate}
\end{theorem}
\emph{Proof}: (i) By results established above, $\CH \cong \LEX \preceq \ACK$.\medskip

\noindent(ii) $\ACK$ differs essentially from $\CH$, $\VN$, and $\Z$ insofar as its terms grow very slowly in rank, so that the lengths of $\ACK$-numbers are very large in comparison to their ranks.  This is the idea underlying the proof we shall now present.  Throughout, let $\N$ be $\VN$, $\Z$, or $\CH$.  The same proof will work for all.

Now suppose, by way of contradiction, that $\varphi$ is a global function that measures $\ACK$ in $\N$: that is,
\[
EA_0\vdash L\mbox{ in }\ACK\, .\supset.\, \varphi(L)\mbox{ in }\N\ \&\  \varphi(L)\simeq_\mathrm{c} L
\]
But by the Rank Bounding Lemma, there is a classical natural number $\mathbf{K}$ such that
\[
EA_0\vdash L\mbox{ in }\ACK\supset\Rank(\varphi(L))<_\mathrm{c} \Rank(L)+_\mathrm{c}\mathbf{K}.
\] 
Now, for any classical natural number $\mathbf{k}$,
\begin{enumerate}
\item[(a)] $\Rank([\cdots, \mathbf{2_k}_{\ACK}]_\ACK)\simeq_\mathrm{c} \mathbf{k}+_\mathrm{c}\mathbf{2}$
\item[(b)] $\Rank([\cdots, \mathbf{2_k}_{\N}]_\N)\simeq_\mathrm{c} \mathbf{2}_\mathbf{k}+_\mathrm{c}\mathbf{2}$
\end{enumerate}
Thus,
\[
\mathbf{2_K+2}\simeq_\mathrm{c} \Rank([\cdots, \mathbf{2_K}_{\N}]_{\N})\simeq_\mathrm{c} \Rank(\varphi([\cdots, \mathbf{2_K}_{\ACK}]))<_\mathrm{c}  \mathbf{K+k+2}
\]
We have thus reached a contradiction and (ii) is established.\medskip

\noindent(iii) Suppose, by way of contradiction, that $\VN \preceq \ACK$. Then, since
\[
\ACK \stackrel{\lambda x\, \mathrm{suplog}_{\mathbf{2}}(x)}{\rightarrow} \Z
\]
we have
\[
\VN \stackrel{\lambda x\, \mathrm{suplog}_{\mathbf{2}}(x)}{\rightarrow} \Z
\]
which contradicts Lemma \ref{VNZconstant}.  Similarly for $\Z$.  \medskip

\noindent This completes our proof.  \hfill$\Box$

\section{Weak subsystems of $\ACK$}\label{wksubsysack}

The systems $\VN$, $\Z$, and $\CH$ all have the weakest closure conditions a natural number system can possibly have. For in these systems the arithmetical function $\lambda x\,x+_\mathrm{c}\mathbf{1}$  is maximally powerful, and every natural number system must be closed under this function. But by taking successive $S$-ary numerals starting with these systems we can obtain an unlimited succession of further systems in which addition, multiplication, and the succession of functions in the $\mu_{\mathbf{n}}$ hierarchy are maximally powerful in the corresponding natural number systems.

Exponentiation is maximally powerful in the  natural number system $\ACK$, which is as far as we can go since all arithmetical functions can be represented there.

These facts naturally raise the question: for which arithmetical global functions $\varphi$ do there exist natural number systems in which $\varphi$ is maximally powerful? This is the question to which we now turn our attention.

We need to concentrate our attention on arithmetical global functions that are \emph{regular above some fixed bound} in the sense of the following definition. 

\begin{definition} Let $\varphi$ be an arithmetical global function and $\mathbf{K}$ be a natural number.  Then, we say that \emph{$\varphi$ is regular above $\mathbf{K}$} if
\begin{enumerate}
\item[\emph{(i)}]$\forall x[\mathbf{K} \leq_\mathrm{c}  x \supset x <_\mathrm{c}  \varphi(x)]$;
\item[\emph{(ii)}]$\forall x\forall y[\mathbf{K} \leq_\mathrm{c}  x \leq_\mathrm{c}  y \supset \varphi(x) -_\mathrm{c} x\leq_\mathrm{c}  \mathbf{2}^y - y]$; and
\item[\emph{(iii)}]$\forall x\forall y[\mathbf{K} \leq_\mathrm{c}  x \leq_\mathrm{c}  y \supset \varphi(x) \leq_\mathrm{c}  \varphi(y)]$.
\end{enumerate}
We call a global arithmetical function that is regular above some fixed bound $\mathbf{K}$ a \emph{regular function}.
\end{definition}

Now we shall carry out a construction which, starting with an arithmetical global function $\varphi$ that is regular above some natural number, yields a natural number system $\ACK_\varphi$ whose terms are $\ACK$-terms, whose numbers are orderings of $\ACK$-terms arranged in their Ackermann order, and in which $\varphi$ can be proved  maximally powerful in $EA$.\footnote{Notice that $\LEX$ is a system of this sort in which $\varphi=\lambda x(x+\mathbf{1})$ can be proved  maximally powerful in $EA$.}

The idea behind our construction starts from the observation that the sequence of terms of $\ACK$
\[
l_0\,(=[\,]=\varnothing),\,l_1(=[s_0]),\,l_2\,(=[s_0,s_1]),\,\cdots, \, l_k\,(=[s_0,\cdots,s_{k-1}]),\,\cdots\,\,\,\footnote{Recall that the $\ACK$-terms are themselves linear orderings, namely finite initial segments of the sequence
\[
s_0\,(=\varnothing),\,s_1\,(=\{\varnothing\}),\,\cdots
\]
of pure sets listed in the order determined by their Ackermann codes.}
\]

\noindent can be obtained from the sequence of terms of $\LEX$
\[
\,l_0,\,l_1,\,l_2,\,l_4,\,l_{16},\,\cdots,\,l_{2_k}\,(=[s_0,\cdots,s_{2_k-1}])\,\cdots
\] 
 
\noindent by interpolating all the missing $\ACK$-terms. We will define the sequence of terms for $\ACK_\varphi$ in a similar manner by interpolating $\ACK$-terms into the sequence of $\LEX$-terms, but leaving some out so that, in general, there are fewer interpolated terms preceding a given $\LEX$-term in the sequence of  $\ACK_\varphi$-terms than in the sequence of $\ACK$-terms:  in fact, we will interpolate exactly enough $\ACK$-terms so that, if there are $M$ terms of $\ACK_{\varphi}$ occurring below $l_{2_k}$, then there are $\varphi(M)$ terms of $\ACK_{\varphi}$ occurring below $l_{2_{k+1}}$, the next term of $\LEX$.  In this way, we ensure that $\ACK_{\varphi}$ is closed under $\varphi$.  Indeed, by the Global Function Bounding Lemma, this will also ensure that $\varphi$ is maximally powerful in $\ACK_{\varphi}$.  The rigorous formulation of the construction just described and the arguments just given is the work of the remainder of this section.

Now for the construction of $\ACK_\varphi$. Suppose we are given a regular function, $\varphi$,  and let $\mathbf{K}$ be the smallest classical natural number above which it is regular. Since $\varphi$ is regular, there is a global function $\ACK\to\ACK$ that represents it in $\ACK$. To simplify the notation let us designate this function by `$\varphi$' as well, so that, using this terminology, $\varphi:\ACK\to\ACK$. Finally, let $\mathbf{N}$ be the smallest classical natural number such that
\[
\forall x[ x\mbox{ in }\ACK\ \&\   x\leq_\mathrm{c} \mathbf{K}\,.\supset.\,\varphi(x)<_\mathrm{c} \mathbf{2_N}]
\]
By clause (i) in the definition of regularity, $\mathbf{K} < \mathbf{2_N}$.

Now we must define the natural number system $\ACK_\varphi$. We shall begin by describing the species of $\ACK_\varphi$-terms (which is naturally ordered by their lengths). Once this is done, we will define the successor function, $\sigma_{\ACK_\varphi}$, for $\ACK_\varphi$. The initial term of $\ACK_\varphi$ is the initial term of $\ACK$, namely, $[\,]$.

To begin with we consider the sequence of $\LEX$ terms
\[
L_0,\,L_1,\,L_2,\,L_3,\,\cdots,\,L_k,\,\cdots
\]
 which is the subsequence
\[
l_0,\,l_1,\,l_2,\,l_4\cdots,\,l_{2_k},\,\cdots
\]
of the sequence $l_0,\,l_1,\,l_2,\,l_3,\,\cdots\,$ of all $\ACK$-terms.

Now we define the sequence of $\ACK_\varphi$-terms inductively in stages  by interpolating new $\ACK$-terms into the sequence of $\LEX$-terms. The number of terms generated by the end of Stage $\mathbf{n}$ is $h_\mathbf{n}$:  our aim is define each stage in such a way that $h_{\mathbf{n}} = \varphi^\mathbf{n}(\mathbf{2_N})$.\bigskip

\noindent\textbf{Stage 0:} We start with \emph{all} the $\ACK$-terms preceding  $l_{\mathbf{2_N}}$
\[
\underbrace{l_0,\,l_1,\,l_2,\cdots,\,l_{\mathbf{2_N}-1}}_{\mbox{Stage $\mathbf{0}$}}
\]
There are $\mathbf{2_N}$ of these terms available at Stage \textbf{0} so we set $h_\mathbf{0}=\mathbf{2_N}$. \bigskip

\noindent\textbf{Stage 1:} Now we add the next $\varphi(h_\mathbf{0})-h_\mathbf{0}\,\,\,(=\varphi(\mathbf{2_N})-\mathbf{2_N}\, )$ terms in the sequence of  $\ACK$ starting with $l_{\mathbf{2_N}}$ to obtain
\[
\underbrace{l_0,\,l_1,\,l_2,\cdots,\,1_{\mathbf{2_N}-1}}_{\mbox{Stage $\mathbf{0}$}},\, \, \,\underbrace{l_{\mathbf{2_N}},\,l_{\mathbf{2_N}+1},\,\cdots,\,l_{\mathbf{2_N}+k_\mathbf{1}-1}}_{\mbox{Stage $\mathbf{1}$}}
\]
Where $k_\mathbf{1}=\varphi(h_\mathbf{0})-h_\mathbf{0}\,\,\, (=\varphi(\mathbf{2_N})-\mathbf{2_N}\, )$. Thus the number of terms  available at the end of Stage \textbf{1} is $h_\mathbf{1}=h_\mathbf{0}+ k_\mathbf{1}\,\,\,(=\varphi(\mathbf{2_N}))$.\bigskip

\noindent\textbf{Stage 2:} Now we skip over the next $2_\mathbf{N+1}-\varphi(\mathbf{2_N})$ terms of $\ACK$, and add $\varphi(h_\mathbf{1})-h_\mathbf{1}\, \, \, (=\varphi^{\mathbf{2}}(\mathbf{2_N}) - \varphi(\mathbf{2_N})$) new terms of $\ACK$, starting with $l_{2_\mathbf{N+1}}$.  Thus by the end of Stage \textbf{2} we have generated the initial segment 
\[
\underbrace{l_0,\cdots,l_{\mathbf{2_N}-1}}_{\mbox{Stage $\mathbf{0}$}},\,\,\,\underbrace{l_{\mathbf{2_N}},\,l_{\mathbf{2_N}+1}\cdots,l_{\mathbf{2_N}+k_\mathbf{1} - \mathbf{1}}}_{\mbox{Stage $\mathbf{1}$}},\,\,\,\underbrace{l_{\mathbf{2_{N+1}}},l_{\mathbf{2_{N+1}+1}},\cdots,l_{\mathbf{2_{N+1}}+k_\mathbf{2} - \mathbf{1}}}_{\mbox{Stage $\mathbf{2}$}}
\]
where $k_\mathbf{2}=\varphi(h_\mathbf{1})-h_\mathbf{1}\,\,\, (=\varphi^{\mathbf{2}}(\mathbf{2_N})-\varphi(\mathbf{2_N})\, )$. The number of terms generated by the end of Stage \textbf{2} is thus $h_\mathbf{2}=h_\mathbf{1}+k_\mathbf{2}\,\,\,\,(=\varphi^\mathbf{2}(\mathbf{2_N}))$.
\[
\vdots \hspace{20mm} \vdots \hspace{20mm} \vdots \hspace{20mm} \vdots \hspace{20mm} \vdots \hspace{20mm} \vdots
\]
\textbf{Stage n+1:} By the end of the Stage $\mathbf{n}$, we have generated the initial segment
\[
\underbrace{l_0,\cdots,l_\mathbf{2_N-1}}_{\mbox{Stage $\mathbf{0}$}},\, \, \, \underbrace{l_\mathbf{2_N}\cdots,l_{\mathbf{2_N} + k_{\mathbf{1}} - \mathbf{1}}}_{\mbox{Stage $\mathbf{1}$}},\cdots, \underbrace{l_\mathbf{2_{N+(n-1)}},l_{\mathbf{2_{N+(n-1)}}+1}\cdots,l_{\mathbf{2_{N+(n-1)}}+k_\mathbf{n}-1}}_{\mbox{Stage $\mathbf{n}$}}
\] 
where $k_{\mathbf{i}} = \varphi(h_{\mathbf{i-1}}) - h_{\mathbf{i-1}} = \varphi^{\mathbf{i}}(\mathbf{2_N}) -  \varphi^{\mathbf{i-1}}(\mathbf{2_N})$.  This initial segment contains $h_\mathbf{n}=\varphi^\mathbf{n}(\mathbf{2_N})$ terms. Now we skip over the next $2_\mathbf{N+n}-\varphi(2_\mathbf{N+(n-1)})$ terms of $\ACK$, and add $\varphi(\mathbf{2_{N+n}})-\varphi(\mathbf{2_{N+(n-1)}})$ new successive $\ACK$-terms, starting with the term $l_{2_\mathbf{N+n}}$ to obtain
\[
\underbrace{l_0,\cdots,l_\mathbf{2_N-1}}_{\mbox{Stage $\mathbf{0}$}}, \cdots, \underbrace{l_\mathbf{2_{N+(n-1)}},\cdots,l_{\mathbf{2_{N+(n-1)}}+k_\mathbf{n}-1}}_{\mbox{Stage $\mathbf{n}$}}, \,\,\, \underbrace{l_\mathbf{2_{N+n}},\cdots,l_{\mathbf{2_{N+n}}+k_\mathbf{n+1}-1}}_{\mbox{Stage $\mathbf{n+1}$}}
\] 
where $k_\mathbf{n+1}=\varphi^\mathbf{n+1}(\mathbf{2_N})-\varphi^\mathbf{n}(\mathbf{2_N})$. The number of terms generated by the end of Stage \textbf{2} is thus $h_\mathbf{n+1}=\varphi^\mathbf{n+1}(\mathbf{2_N})$.

It is guaranteed that there are at least $k_{\mathbf{n+1}}$ terms of $\ACK_\varphi$ in the sequence $l_{\mathbf{2_{N+n}}}, \cdots, l_{\mathbf{2_{N+n+1} - 1}}$, since this sequence contains $\mathbf{2_{N+n+1}} -\mathbf{2_{N+n}}  = \mathbf{2^{2_{N+n}}} -\mathbf{2_{N+n}}$ terms of $\ACK$ and, since $\varphi$ is regular above $\mathbf{K}$, 
\[
k_\mathbf{n+1}=\varphi^\mathbf{n+1}(\mathbf{2_N})-\varphi^\mathbf{n}(\mathbf{2_N}) <_\mathrm{c}  \mathbf{2^{2_{N+n}}} -\mathbf{2_{N+n}}
\]
We can continue in this manner to generate longer and longer initial segments of the sequence of $\ACK_\varphi$ terms in their natural  $\ACK$ order. Of course this (informal) inductive definition doesn't, as it stands, define the natural number system $\ACK\varphi$. To give a proper definition of $\ACK\varphi$ we must give a direct definition of its successor function $\sigma_{\ACK_\varphi}$.

Our strategy will be to define a function, $\Stage_\varphi$, such that, given a $\CH$-number $[V_0, \cdots, V_{\mathbf{N}}, \cdots, V_{\mathbf{N} + n}]$,
\begin{eqnarray*}
\Stage_{\varphi}(V_0) =  \cdots = \Stage_\varphi(V_\mathbf{N}) & =&  [l_0, \cdots, l_{\mathbf{2_N - 1}}] \\
\Stage_{\varphi}(V_\mathbf{N + 1}) & = & [l_{\mathbf{2_N}},\,l_{\mathbf{2_N}+1},\,\cdots,\,l_{\mathbf{2_N}+k_\mathbf{1}-1}] \\
\Stage_{\varphi}(V_{\mathbf{N + 2}}) & = & [l_{\mathbf{2_{N+1}}},l_{\mathbf{2_{N+1}+1}},\cdots,l_{\mathbf{2_{N+1}}+k_\mathbf{2} - \mathbf{1}}] \\
\vdots & & \vdots \\
\Stage_{\varphi}(V_{\mathbf{N} + n}) & = & [l_{\mathbf{2}_{\mathbf{N} + n-1}}, l_{\mathbf{2}_{\mathbf{N} + n-1}+\mathbf{1}}\cdots,l_{\mathbf{2}_{\mathbf{N} + n-1}+k_n -\mathbf{1}}]
\end{eqnarray*}
Here $\Stage_\varphi(V_{\mathbf{N + k}})$ corresponds to the sequence of terms added at Stage $\mathbf{k}$ in the informal construction above.

We already know that, from $[V_0, \cdots, V_{\mathbf{N}}, \cdots, V_{\mathbf{N} + n}]$, we can obtain the full $\ACK$ term,
\[
[l_0, \cdots, l_{\mathbf{2_N - 1}}, l_{\mathbf{2_N}},\,l_{\mathbf{2_N}+1},\,\cdots\, l_{\mathbf{2}_{\mathbf{N} + n}}]
\]
Thus, the only hurdle to internalizing our definition of $\Stage_\varphi$ is internalizing the definition of $k_n$.

Following our informal construction above, we want to define $k_n$ such that $k_n = \varphi^n(\mathbf{2_N}) - \varphi^{n-1}(\mathbf{2_N})$.  We do this by first defining the local function $\lambda n\, \varphi^n(\mathbf{2_N}) : [0, \cdots,  k, k+1, \cdots,n] \rightarrow \mathrm{Card}(V_{\mathbf{N} + n + 1})$ by recursion along the linear ordering $[0, \cdots, k, k+1, \cdots, n]$, which is given to us in the $\mathcal{CH}$-number $[V_0, ..., V_\mathbf{N}, ..., V_{\mathbf{N} + n}$]:
\begin{eqnarray*}
\varphi^0(\mathbf{2_N}) & = & \mbox{the unique $x$ in $\mathrm{Card}(V_{\mathbf{N} + n + 1})$ such that $x \simeq_\mathrm{c} \mathbf{2_N}$} \\
\varphi^{k+1}(\mathbf{2_N}) & = & \mbox{the unique $x$ in $\mathrm{Card}(V_{\mathbf{N} + n + 1})$ such that $x \simeq_\mathrm{c} \varphi(\varphi^{k}(\mathbf{2_N}))$}
\end{eqnarray*}
The efficacy of this definition is guaranteed, for, by clause (ii) in the definition of regularity, if $a \leq_\mathrm{c}  V_{\mathbf{N}+k}$, then
\[
\varphi(a) \leq_\mathrm{c} \mathbf{2}^{V_{\mathbf{N}+k}} \simeq_\mathrm{c}  V_{\mathbf{N}+k+1} \subseteq V_{\mathbf{N} + n + 1}
\]
Thus, given $[V_0, \cdots, V_{\mathbf{N}}, \cdots, V_{\mathbf{N} + n}]$, we may define
\[
k_n =_{\mathrm{df.}} \varphi^n(\mathbf{2_N}) - \varphi^{n-1}(\mathbf{2_N})
\]
which completes our definition of $\Stage_\varphi$.

Thus, we have a global function that takes $[V_0, \cdots, V_{\mathbf{N}}, \cdots, V_{\mathbf{N}+n}]$ in $\CH$ first to $k_n$ and then to
\[
[\Stage_\varphi(V_0), \cdots, \Stage_\varphi(V_\mathbf{N}), \cdots, \Stage_\varphi(V_{\mathbf{N}+n})].
\]
Thus, let
\[
\mbox{\textsc{Stage}}_\varphi(V_{\mathbf{N} + n}) =_{\mathrm{df.}} \Stage_\varphi(V_\mathbf{N}) * \cdots * \Stage_\varphi(V_{\mathbf{N}+n}) 
\]
With this in hand, we can give a rigorous definition of $\ACK_\varphi$: 
\begin{definition}
Let $\ACK_\varphi$ be the iteration system generated from $l_0 = [\, ]$ by the following iteration function:  Given an $\ACK$-term, $l_k$, if $L_{\mathbf{N} + n+1}$ is a term of $\LEX$ such that $l_k \subsetneqq_* L_{\mathbf{N} + n+1}$, and, if $l_k$ is in the field of the concatenation
\[
\mbox{\textsc{Stage}}_\varphi(V_{\mathbf{N} + n+\mathbf{2}}) = \Stage_\varphi(V_\mathbf{N}) * \cdots * \Stage_\varphi(V_{\mathbf{N} + n + 1}) * \Stage_\varphi(V_{\mathbf{N} + n+\mathbf{2}})
\]
then let
\[
\sigma_{\ACK_\varphi}(l_k) =_{\mathrm{df.}} \Next_{\mbox{\textsc{Stage}}_\varphi(V_{\mathbf{N} + n+\mathbf{2}})}(l_k)
\]
\end{definition}
Note that it follows from this definition that, for any term $V_{\mathbf{N}+n}$ of $\CH$, $\mbox{\textsc{Stage}}_\varphi(V_{\mathbf{N}+n})$ is a number of $\ACK$.  Moreover,
\begin{enumerate}
\item[(a)] $\mbox{\textsc{Stage}}_\varphi(V_\mathbf{N}) =_{\mathrm{df.}} \Stage_\varphi(V_{\mathbf{N}}) \simeq_\mathrm{c} \mathbf{2_N}$;
\item[(b)] $\mbox{\textsc{Stage}}_\varphi(V_{\mathbf{N}+n+\mathbf{1}}) =_{\mathrm{df.}} $

\hspace{20mm} $\mbox{\textsc{Stage}}_\varphi(V_{\mathbf{N}+n}) * \Stage_\varphi(V_{\mathbf{N}+n+\mathbf{1}}) \simeq_\mathrm{c} \varphi(\mbox{\textsc{Stage}}_\varphi(V_{\mathbf{N}+n}))$;
\item[(c)] $\Rank(\mbox{\textsc{Stage}}_\varphi(V_{\mathbf{N}+n+1})) \simeq_\mathrm{c} \Rank(L_{\mathbf{N}+n+1})$, $\Rank(\mbox{\textsc{Stage}}_\varphi(V_{\mathbf{N}+n})) \simeq_\mathrm{c} \Rank(L_{k})$, and $\Rank(L_{\mathbf{N}+n+1}) \simeq_\mathrm{c} \Rank(L_{\mathbf{N}+n}) + \mathbf{1}$, so
\[
\Rank(\mbox{\textsc{Stage}}_\varphi(V_{\mathbf{N}+n+1})) \simeq_\mathrm{c} \Rank(\mbox{\textsc{Stage}}_\varphi(V_{\mathbf{N}+n})) + \mathbf{1}
\]
\end{enumerate}

\begin{theorem}Suppose $\varphi$ is regular above $\mathbf{K}$. Then $\varphi$ is maximally powerful in $\ACK_\varphi$.
\end{theorem}

\noindent\emph{Proof}: First we must  show that $\ACK_\varphi$ is closed under $\varphi$. For this it will suffice to exhibit a global function $\gamma$ such that
\[
EA_0\vdash L\mbox{ in }\ACK_\varphi\, .\supset.\, \gamma(L) \mbox{ in } \ACK_\varphi \ \&\  \varphi(L)\leq_\mathrm{c} \gamma(L)
\]
To that end, define $\gamma$ as follows.  If $L$ is a number of $\ACK_\varphi$, then there are two cases:
\begin{enumerate}
\item[(1)] If $L \leq_\mathrm{c} \mathbf{K}$, then let
\[
\gamma(L) =_\mathrm{df.} \mbox{\textsc{Stage}}_\varphi(V_\mathbf{N})
\]
\item[(2)] If $\mathbf{K} <_\mathrm{c} L$, then we can recover the $\mathcal{CH}$-number, $[V_0, ..., V_\mathbf{N}, ..., V_{\mathbf{N} + k}]$ such that $\Field(L) \subseteq V_{\mathbf{N} + k}$.  Thus, $L \subseteq_*  \mbox{\textsc{Stage}}_\varphi(V_{\mathbf{N}+k})$.  In this case, let
\[
\gamma(L) =_{\mathrm{df.}} \mbox{\textsc{Stage}}_\varphi(V_{\mathbf{N}+k+\mathbf{1}})
\]
\end{enumerate}
If case (1) obtains, then by (a) above
\[
\varphi(L) <_\mathrm{c}  \mathbf{2_N} \simeq_\mathrm{c} \mbox{\textsc{Stage}}_\varphi(V_\mathbf{N}) =_{\mathrm{df.}} \gamma(L)
\]
If case (2) obtains, then, by clause (ii) in the definition of regularity and by (b) above,
\[
\varphi(L) \leq_\mathrm{c}  \varphi(\mbox{\textsc{Stage}}_\varphi(V_{\mathbf{N}+k})) \simeq_\mathrm{c} \mbox{\textsc{Stage}}_\varphi(V_{\mathbf{N}+k+\mathbf{1}}) =_{\mathrm{df.}} \gamma(L)
\]
Thus, $\ACK_\varphi$ is closed under $\varphi$. \medskip

Next, we show that, if $\ACK_\varphi$ is closed under an arithmetical global function $\psi$, there is a classical natural number $\mathbf{n}$ such that $\psi \ll \varphi^{\mathbf{n}}$.  As usual, we use `$\psi$' to name the function that represents $\psi$ in $\ACK_\varphi$.  Then, by the Rank Bounding Lemma, there is a classical natural number $\mathbf{n}$ such that
\[
EA_0 \vdash \Rank(\psi(a)) \leq_\mathrm{c}  \Rank(a) + \mathbf{n}
\]
Now, suppose $L$ is a number of $\ACK_\varphi$, suppose $\mathbf{K} < L$, and suppose $\mbox{\textsc{Stage}}_\varphi(V_{\mathbf{N}+k}) \subseteq_* L \subseteq_* \mbox{\textsc{Stage}}_\varphi(V_{\mathbf{N}+k+1})$.  Then
\[
\Rank(L) \simeq_\mathrm{c} \Rank(\mbox{\textsc{Stage}}_\varphi(V_{\mathbf{N}+k+1})).
\]
  By (c) above,
\[
\Rank(\psi(L)) \leq_\mathrm{c}  \Rank(L) + \mathbf{n} \leq_\mathrm{c} \Rank(\mbox{\textsc{Stage}}_\varphi(V_{\mathbf{N}+k+ \mathbf{1} + \mathbf{n}}))
\]
which entails
\[
\psi(L) \subseteq_* \mbox{\textsc{Stage}}_\varphi(V_{\mathbf{N}+k+\mathbf{n}}) \simeq_\mathrm{c} \varphi^{\mathbf{n}}(\mbox{\textsc{Stage}}_\varphi(V_{\mathbf{N}+k})) \leq_\mathrm{c}  \varphi^{\mathbf{n}}(L)
\]
since $\varphi$ is regular above $\mathbf{K}$ and $\mbox{\textsc{Stage}}_\varphi(V_{\mathbf{N}+k}) \subseteq_* L$. \hfill $\Box$

{\footnotesize

\bigskip

\end{document}